\documentclass[twocolumn]{autart}    

\usepackage{graphicx}          

\usepackage[english]{babel}
\usepackage{amsmath}
\usepackage{amssymb}
\usepackage{amsfonts}
\usepackage{mathtools} 
\usepackage{MnSymbol}

\usepackage{algorithm}
\usepackage{algpseudocode}
\usepackage{textcomp}
\usepackage{color}

\newtheorem {defi}{Definition}
\newtheorem {theorem}{Theorem}
\newtheorem {remark}{Remark}
\newtheorem {lemma}{Lemma}
\newtheorem {example}{Example}

\begin{document}

\begin{frontmatter}

\title{
On the distributed backward reachability problem for large scale systems 
\thanksref{footnoteinfo}} 

\thanks[footnoteinfo]{This paper was not presented at any IFAC 
meeting. Corresponding author D.~Liuzza. Tel. +39 06 94005727. }

\author[Enea]{Davide Liuzza}\ead{davide.liuzza@enea.it},    
\author[Unimore,Chalmers]{Paolo Falcone}\ead{falcone@unimore.it},               
\author[Unisannio]{Massimo Tipaldi}\ead{mtipaldi@unisannio.it},
\author[Unisannio]{Luigi Glielmo}\ead{glielmo@unisannio.it}  

\address[Enea]{ENEA, Fusion and Nuclear Safety Department, Frascati, 00044 Italy}  
\address[Unimore]{University of Modena and Reggio Emilia,
Engineering Department “Enzo Ferrari”,
Modena, Italy}             
\address[Chalmers]{Chalmers University of Technology,
Department of Electrical Engineering,
G\"oteborg, Sweden}           
\address[Unisannio]{
University of Sannio, Department of Engineering, Piazza Roma, 82100 Benevento, Italy}

\begin{keyword}                           
distributed reachability, distributed  control, networked control systems             
\end{keyword}                             

\begin{abstract}                          
Backward reachability (also termed controllability) has been extensively studied  in control theory, and tools for a wide class of systems have been developed. Nevertheless, assessing a backward reachability analysis or synthesis remains  
challenging as the system dimension grows.
In this paper we study the backward reachability problem for large scale networked nonlinear systems 
with coupled dynamics and subject to states and inputs nonlinear constraints.
We propose a theory for completely general nonlinear  constrained large scale controllability problems. 
We demonstrate that it is always possible to recast such problems for the overall large scale system into an equivalent distributed form where, without 
introducing any conservativeness, each node of the network iteratively solves a local reachability subproblem  
by exchanging information with the adjacent  
nodes.  Although the proposed algorithm is completely decentralized, the solution of the backward reachability problem for the overall system is equivalently determined by the local ones and satisfies all the given constraints.
Not being linked to any specific assumption on the system dynamics nor static constraints, the proposed results hold irrespectively of any possible analytical/numerical solver to be adopted for backward reachability computation. 
\end{abstract}

\end{frontmatter}

\section{Introduction}
\subsection{Motivation and relevant literature}\label{subsec:motivation&literature}
{\em Backward reachability} analysis for dynamical systems, also termed as {\em controllability}, has been widely studied by control theoreticians. Indeed, knowing how to reach a desired state space region of a dynamical system from a starting region through an admissible control action is a powerful tool enabling a number of theoretical analyses and practical applications. 
The reader interested in fundamental results about reachability analysis is referred to
\cite{demi:69,glsc:71,berh:71,be:72,begi:00,miba:05,rake:06,abam:07}.
In particular, reachability has been studied and exploited as a fundamental tool for evaluating or enforcing state space invariance of control systems  \cite{berh:71,la:93,bl:99,rake:05,maka:13} or reach-avoid set control  and differential games \cite{miba:05,dito:10,maly:11,bozi:11,fevi:20}. A renewed interest in the (both forward and backward)  reachability problem is witnessed by more recent literature too, where this tool is exploited in the derivation of dynamical systems abstraction techniques and symbolic control approaches for the verification of fundamental properties such as safety or for the enforcement of formal logics specifications \cite{tapa:06,pogi:08,woto:09,yobe:10,zapo:11,yotu:11,nioz:12,wotoIROS:13,cogo:15,bodi:17,medi:18,medi:19}.
This fact, indeed, comes from its inherent peculiarity of addressing how two regions of the state space (a starting and an ending region) are mapped through the dynamics of a (in general nonlinear) system under selected inputs.
In this context, forward reachability addresses the problem of propagating a starting state space region through the dynamics of the system under given inputs (or in free evolution), so as to obtain the system's dynamical flow. Backward reachability addresses instead the opposite problem of finding the starting region and the control inputs so as to reach, under the system dynamical flow, a target ending region.
\newline
These two different aspects of the same concept find, in general, different applications since while forward reachability is mostly useful for simulation, safety and verification applications (the latter can be, however, covered also via the backward one), backward reachability is particularly useful for control synthesis applications, e.g., for setting the starting, terminal and input constraints of a Model Predictive Control (MPC) problem or for invariant set computation, among others.

A modern trend in automatic control is to address the analysis and control of complex, high dimensional and interconnected dynamical systems. 
Examples of such problems are abundant in power networks, biological systems, power management, transportation problems and robotics \cite{amwo:05,ghto:05,yobe:10,ozto:11,cogo:15,bodi:17}. In this regard and for the reasons reported above, reachability analysis appears to be a promising candidate tool to assess formal properties of these complex systems, derive related abstraction techniques and design control laws. 
 
In general, however, several existing reachability theoretical or numerical tools developed in the classical control do not scale up properly
for those complex systems.
\newline
In particular for the backward reachability case, 
a general theory is missing that is able to convert a large scale system controllability problem as a whole into an equivalent distributed one. 
\newline
This point is of particular importance since, when looking at a large scale system as a unique high-dimensional one, existing algorithms lead to intractable numerical problems.


To highlight the paper contribution and to provide a (non exhaustive) introduction to the context, in what follows we give a brief overview on the reachability literature.
Specifically, when coming to forward reachability, recent numerical tools coping with some classes of nonlinear systems  and adopting set approximation and over-approximation techniques have been developed  \cite{dumi:15,koga:15,chab:13,al:15,scal:20}. Although several of these tools do not explicitly focus on the scalability problem, they provide rather good performances in assessing the reachable set of relatively high dimensional systems.

Also, a good scalability for piecewise affine linear systems over template polytopes is shown in \cite{frle:11}.

Reachability problem for linear systems are also addressed in \cite{kuva:00} via inner ellipsoidal approximations of reachable tubes. Ellipsoidal reachability for linear systems is also addressed in \cite{babi:20}, both for state estimation and control.

Relevant numerical codes for (forward) reachability are reported at the link  \cite{CSShybrid}. In particular some of them, such as SpaceEx and MATISSE for linear systems of a certain specific structure and JuliaReach for nonlinear systems, show good scalability properties at the system dimension growth. To achieve fast computation, approximations of the considered sets or the system dynamics are considered. In addition to this, a suitable code structure allows to leverage on accelerated GPU support.

Reachability and control invariance have been addressed, among others, in \cite{maka:13} where theoretical results are integrated into already existing reachability tools to provide reachability algorithms able to operate with systems degree up to some tens.  

Backward reachability for certain classes of nonlinear system is addressed in \cite{mava:14} via a semidefinite programming approach in order to simultaneously compute both the backward reachable set and the appropriate control action.

Distributed invariant set computation through reachability for linear affine systems 
with coupled dynamics is addressed in \cite{nioz:16}, where separable invariant sets are computed for each subsystem of the network, based on an assume-guarantee approach. However, this approach only applies to linear affine systems with state coupling.
 
Hamilton-Jacobi formulation along with system decomposition has been successfully applied in \cite{chhe:16} to solve the backward reachability for high-dimensional nonlinear system. 
Despite the significant improvement of the overall computational time, the adopted approach is inherently conservative since the couplings in the dynamics are treated as disturbances.
 
Other papers have been recently proposed trying to address scalability issues under specific problem structures \cite{mi:11,coar:15}. An early, interesting Hamilton-Jacobi formulation of the reachability problem through projection can be found in \cite{mito:03}, where the reachable set of a high-order system is computed through projections over lower dimensional spaces. Each projection results in an over-approximation of the projection of the exact reachable set and, similarly, the intersection of the back-projection of the formers is an over-approximation of the original overall reachable set. 

All the cited approaches, however, present limitations both/either in generalizing to arbitrarily large system and general nonlinear system structures and/or in the conservativeness of the resulting reachable set. 
 
An alternative projection approach has recently been proposed  in \cite{chhe:18}, where backward reachability is investigated for a special class of nonlinear systems, which admit a state space decomposition such that each resulting subsystem evolution depends only on a subset of the overall state variables and on a certain number of common state variables shared among all of the other subsystems. 
Backward reachable sets are then locally computed for each subsystem and then reconstructed for the overall nonlinear system through a centralized union/intersection operation.

\subsection{Paper contribution} \label{subsec:contribution}
As it is possible to recognize from the above referenced papers, different approaches have been adopted for the reachability problem, spanning from theoretical to numerical investigations.
\newline 
When coming to large scale systems, several papers approach the problem from a numerical viewpoint, looking for set and dynamics approximations or for algorithms parallelization to speed up the computations. 
\newline
In our paper, we study the reachability problem from a completely different viewpoint and resort to a newly developed theoretical approach. 
\newline
Specifically, we address the backward reachability problem for large scale systems arising by the interconnection of several dynamical {\em agents} possibly linked via their dynamics (both states and inputs) and/or via static constraints. 
For these kind of systems, and totally generic in terms of dynamics and constraints, we study whether a backward reachability problem given for the overall large scale system can be decomposed into local subproblems.
\newline
We provide a theory and a methodology showing that this is always possible. In particular, we prove that for any generic controllability problem it is always possible to find a distributed representation via suitably exploiting the patterns already present in the overall system and derive a distributed reformulation of the problem. 
\newline
Such distributed reformulation is equivalent to the one of the overall system (that here we will term {\em centralized}) in the sense that, combining the solutions of the local reachability problems at each agent level it is possible to reconstruct the centralized solution with no information losses. 
Nevertheless, although doable, such reconstruction does not need to be performed. Indeed, according to the proposed decomposition, each agent will only solve the minimal ``portion'' of the centralized controllability problem that allows, at the same time, not losing any information with respect to the centralized problem and not adding any further information to the local problem that is not essential for the agent. 
\newline
From a broader viewpoint, the paper investigates the centralized vs distributed control relation. Such a topic has been hugely studied by several authors and represents an important (both theoretical and practical) scientific research line. 
Via the abstract and purely theoretical approach, we show that a general control problem in abstract terms (reaching a target region from a starting one under a control) for large scale systems can always be decomposed into a distributed form. This holds for general nonlinear dynamics and general static constraints. To the best of our knowledge, such general result is not present in the current literature. 
\newline
Therefore, differently from other works in the literature, we start from a theoretical problem and related research question of investigating the centralized vs distributed relation for control and reach a theoretical scheme that is able to provide a distributed solution for any  controllability problem. 
\newline
In view of such an approach, our paper differs from the valuable works presented above, where authors were mainly interested in speeding up parallel/approximated calculation for (forward) reachability. 
Instead, we are interested in a theoretical problem and do not focus on implementation aspects. Specifically, in the paper pseudo-algorithms are provided that convert an abstract backward reachability problem into distributed reachability problems of reduced size (with a number of decision variable that only depends on each local agent neighbourhood).
Each agent is assumed to have communication and computation capabilities, to compute such low-order backward reachable sets and to perform local set operations (union and intersection). The results of the local computations are then iteratively exchanged with the neighbouring agents until obtaining, for each agent, the local reachable set which, together with the ones computed by the other agents, represents an equivalent fully decentralized version of the overall system's reachable set obtained through a centralized computation.
\newline
We wish to emphasize that, in line with the theoretical approach  developed,  we coherently do not address how the local reachability problems are solved. Indeed, we do not focus on some specific system dynamics and, therefore, we cannot provide a specific solver for the proposed reachability. Nevertheless, this point does not represent a criticality since the local reachability can be addressed, if possible, with any of the compatible approaches already developed in the literature and where already results are provided. 
In this sense, the paper contribution is not about proposing some reachability solver (such research line is already well investigated in the literature) but rather providing the theory that allows to distribute large scale problems. In this sense, if the distributed local reachability problems can be solved with some existing method (analytically or numerically), then the large scale system (for any dimension size) can be solved as well. 
\newline
To ease the paper comprehension, a discussion section is provided where ideas on how to cast the proposed results to some specific problem setting are discussed. In this regard, we illustrate with more details the linear affine system case, where no results are previously available in the literature for a completely general setting.  
Also, a numerical example illustrating the main proposed algorithm is provided to allow the reader grasping the theoretical concepts.

Compared to the existing literature, our work contributes to advance the state of the art by:
\begin{itemize}
\item developing a general theory without restrictions to  special classes of nonlinear systems;
\item showing that any controllability problem for large scale systems admits a distributed formulation;
\item fully decentralizing the calculation of the backward reachable set, thus extending the tool to an arbitrarily large number of agents of a networked system (in other words an arbitrarily high-order nonlinear system);
\item allowing state and input nonlinear constraints;
\item avoiding any centralized operations, such as in \cite{mito:03,chhe:18};
\item not introducing any conservativeness. 
\end{itemize}

As already said, and in line with some other papers in the literature, such as \cite{maka:13,chhe:18}, this work solely focuses on the (quite involved) distributed formulation of the backward reachability problem for a constrained networked system. As such, the paper content is mainly devoted to the presentation of the theoretical foundation body, while its applications to specific classes of systems will be presented in dedicated papers.

The organization of the manuscript is as follows. 
In Section \ref{sec:extrusion_generated_set} we provide important preliminary mathematical concepts and related notation that will be heavily exploited in the rest of the paper. The dynamical system model, backward reachability concepts and the problem statement are provided in Section \ref{sec:problem_statement}.     
In Section \ref{sec:distributed_extrusion_generated_set} we move again back on the mathematical formalism provided in Section \ref{sec:extrusion_generated_set} and further develop background results on distributed computation. The latter will be then applied to large scale dynamical systems in Section 
\ref{sec:distributed_reachability}.
The alternate scheme Section \ref{sec:extrusion_generated_set} and Section \ref{sec:distributed_extrusion_generated_set} for mathematical background concepts and Section 
\ref{sec:problem_statement} and Section \ref{sec:distributed_reachability} for the dynamical systems case is adopted to timely provide in Section  \ref{sec:problem_statement} the formal problem statement. In this regard, Section \ref{sec:extrusion_generated_set} only reports those mathematical concepts that are strictly needed for this. Further additional background, essential for the problem resolution but not for the problem statement, is therefore moved later in Section \ref{sec:distributed_extrusion_generated_set}. 
\newline
To complete the paper, in Section \ref{sec:discussion} we discuss the applicability of the theory we developed, by providing some insight into how it can be used for further researches. In particular, we propose ideas on the linear affine system case, sufficiently rich to illustrate the developed theory. 
Conclusions and final remarks are given  in Section \ref{sec:conclusions}.
\newline
Examples are reported to illustrate the concepts and the operators introduced along the paper. 
\newline
To ease the paper reading, examples are reported in Appendix A, while proofs for theorems and lemmas are given in Appendix B.

\section{Extrusion generated set}
\label{sec:extrusion_generated_set}

In preparation of the formal statement of the reachability problem in Section \ref{sec:problem_statement}, in this section 
we preliminary give the concept of {\em extrusion generated set} and
related working operators. 
The concepts here presented will be not only useful for defining the problem addressed in this manuscript, but will be also instrumental to its resolution.

\begin{defi}\label{def:axis_set}
A finite set $\mathcal{B}\subset\Bbb{N}^+$ is said to be an {\em axis set} if it can be written as $\mathcal{B}=\{\beta_1,\beta_2,\dots,\beta_q\}$ such that $\beta_i < \beta_j$ iff $i<j$.
\end{defi}

We provide next the definitions of two useful operators heavily exploited in the rest of the paper. 
\begin{defi}\label{def:projection}
Let us consider two axis sets $\mathcal{B}_1=\{\beta_1^{1},\dots,\beta_{q_1}^1\}$ and $\mathcal{B}_2=\{\beta_1^{2},\dots,\beta_{q_2}^2\}$, with $\mathcal{B}_1\subseteq \mathcal{B}_2$. The {\em projection operator} 
$\mathcal{P}_{\mathcal{B}_1}^{\mathcal{B}_2}(v)$, with $v\in\Bbb{R}^{|\mathcal{B}_2|}$, is defined as
\begin{equation}\label{eq:projection}
\mathcal{P}_{\mathcal{B}_1}^{\mathcal{B}_2}(v):=\{w\},
\end{equation}
with $w=(w^{(1)},\dots,w^{(q_1)})^T\in\Bbb{R}^{|\mathcal{B}_1|}$ such that   $w^{(i)}=v^{(j)}$ iff $\beta^{1}_i=\beta^{2}_j$.

The above definition can be extended (with a slight abuse of notation) to the case where the argument of operator $\mathcal{P}_{\mathcal{B}_1}^{\mathcal{B}_2}(\cdot)$ is a set $V\subseteq \Bbb{R}^{|\mathcal{B}_2|}$. In such a case, we define
\begin{equation*}
\mathcal{P}_{\mathcal{B}_1}^{\mathcal{B}_2}(V):=\bigcup_{v\in V}\mathcal{P}_{\mathcal{B}_1}^{\mathcal{B}_2}(v).
\end{equation*} 
Finally, we impose the following conventions: $\mathcal{P}_{\mathcal{B}_1}^{\mathcal{B}_2}(\emptyset)=\mathcal{P}_{\emptyset}^{\mathcal{B}_2}(V)=\mathcal{P}_{\emptyset}^{\emptyset}(\emptyset)=\emptyset$, $V\subseteq \Bbb{R}^{|\mathcal{B}_2|}$.
\end{defi}

Note that, when projecting a single point in $\Bbb{R}^{|\mathcal{B}_2|}$,  the result of the operator \eqref{eq:projection} is a singleton set. With a slight abuse of notation, in the rest of the paper we will either refer to such set or to the element it contains. It will be clear from the context which of the two meanings we refer to. 
Example \ref{exa:projection_operator} clarifies how the projection operator is computed.

Similarly to Definition \ref{def:projection}, the dual {\em extrusion operator} is defined as 
\begin{defi}\label{def:extrusion}
Let us consider two axis sets $\mathcal{B}_1=\{\beta_1^{1},\dots,\beta_{q_1}^1\}$ and $\mathcal{B}_2=\{\beta_1^{2},\dots,\beta_{q_2}^2\}$, with $\mathcal{B}_1\subseteq \mathcal{B}_2$. The {\em extrusion operator} 
$\mathcal{E}_{\mathcal{B}_1}^{\mathcal{B}_2}(w)$, with $w\in\Bbb{R}^{|\mathcal{B}_1|}$, is defined as
\begin{align}
\mathcal{E}_{\mathcal{B}_1}^{\mathcal{B}_2}(w):=&\left\{v\in\Bbb{R}^{|\mathcal{B}_2|}: v^{(i)}=w^{(j)} \mathrm{\,\, iff \,\,} \beta_i^2=\beta_j^1; \right. \nonumber \\
& \left. v^{(i)}\in\Bbb{R} \mathrm{\,\, otherwise}  \right\}.\label{eq:extrusion}
\end{align}

The above definition can be extended (with a slight abuse of notation) to the case where the argument of the operator $\mathcal{E}_{\mathcal{B}_1}^{\mathcal{B}_2}(\cdot)$ is a set $W\subseteq \Bbb{R}^{|\mathcal{B}_1|}$. In such a case, we define
\begin{equation*}
\mathcal{E}_{\mathcal{B}_1}^{\mathcal{B}_2}(W):=\bigcup_{w\in W}\mathcal{E}_{\mathcal{B}_1}^{\mathcal{B}_2}(w).
\end{equation*} 
Finally, we impose the following conventional definitions: $\mathcal{E}_{\mathcal{B}_1}^{\mathcal{B}_2}(\emptyset)=\mathcal{E}_{\emptyset}^{\emptyset}(\emptyset)=\emptyset$. 
\end{defi}

Similarly to the projection operator, the extrusion operator is illustrated in Example \ref{exa:extrusion_operator}.

Now, let us define 
\begin{equation}\label{eq:B_bar}
\bar{\mathcal{B}}=\bigcup_{i=1}^N \mathcal{B}_i.
\end{equation}

We can provide the following important definition. 

\begin{defi}\label{def:extrusion_generated_set}
Let us consider a collection of $N$ non empty axis sets $\mathcal{B}_1,\mathcal{B}_2,\dots,\mathcal{B}_N$, set $\bar{\mathcal{B}}$ as in \eqref{eq:B_bar}, and a set $\bar{S}\subseteq\Bbb{R}^{|\bar{\mathcal{B}}|}$. The set $\bar{S}$ is called an {\em extrusion generated set} associated to the collection $\mathcal{B}_1,\dots,\mathcal{B}_N$ iff there exist 
$S_i\subseteq \Bbb{R}^{|\mathcal{B}_i|}$, for $i=1,\dots, N$, such that
\begin{equation}\label{eq:extrusion_generated_set}
\bar{S}=\bigcap_{i=1}^N \mathcal{E}_{\mathcal{B}_i}^{\bar{\mathcal{B}}}\left( S_i \right).
\end{equation}
\end{defi}
To provide a clue on the centralized operation performed in \eqref{eq:extrusion_generated_set}, a toy example (namely Example \ref{exa:centralized_computation}) is provided.

According to Definition~\ref{def:extrusion_generated_set}, we denote by
\begin{equation}\label{eq:projection_extrusion_generated_set}
\bar{S}_i=\mathcal{P}_{\mathcal{B}_i}^{\bar{\mathcal{B}}}(\bar{S}),
\end{equation}
the projection of $\bar{S}$ with respect to the axis set $\mathcal{B}_i$. The following theorem holds.

\begin{theorem}\label{thm:bar_S_centralized_from_projection}
Let us consider a collection of $N$ non empty axis sets $\mathcal{B}_1,\mathcal{B}_2,\dots,\mathcal{B}_N$, and an  extrusion generated set $\bar{S}\subseteq\Bbb{R}^{|\bar{\mathcal{B}}|}$ associated to such collection as in Definition~\ref{def:extrusion_generated_set}. Also, let us consider its projection $\bar{S}_i$ as in   \eqref{eq:projection_extrusion_generated_set}. 
Then,
\begin{equation}\label{eq:bar_S_centralized_from_projection}
\bar{S}=\bigcap_{i=1}^N \mathcal{E}_{\mathcal{B}_i}^{\bar{\mathcal{B}}}\left(\bar{S}_i \right). 
\end{equation}

Furthermore, in case $\bar{S}\neq \emptyset$, consider  sets $\hat{S}_i\subseteq \bar{S}_i$, for $i=1,\dots, N$. If there exists at least one
index $j$ and the corresponding set $\hat{S}_j$ such that the strict inclusion holds, i.e., $\hat{S}_j\subset \bar{S}_j$, then $\hat{S}=\bigcap_{i=1}^N \mathcal{E}_{\mathcal{B}_i}^{\bar{\mathcal{B}}}(\hat{S}_i)\subset \bar{S}$.
\end{theorem}

The findings of Theorem \ref{thm:bar_S_centralized_from_projection} are shown in Example \ref{exa:centralized_computation_from_projection}.

Although simple in its form, Theorem \ref{thm:bar_S_centralized_from_projection} plays an important role for the theory we propose later. Its interpretation is given in the following remark.
\begin{remark}\label{rem:bar_S_centralized_from_projection}
From Theorem \ref{thm:bar_S_centralized_from_projection} we have that an extrusion generated set $\bar{S}$ can be obtained from its projections by applying formula  \eqref{eq:bar_S_centralized_from_projection}. Furthermore, the projections $\bar{S}_i$ are the ``tightest'' sets able to generate $\bar{S}$. Any other set $\hat{S}_j\subset \bar{S}_j$, for some $j\in\{1,\dots, N\}$ implies an extrusion generated set $\hat{S}\subset \bar{S}$. 
Such an aspect is particularly important since $\bar{S}$ can be uniquely reconstructed through the distributed information of the $\bar{S}_i$'s stored at the nodes level. Therefore, despite in general a single node in the network is not able to reconstruct alone the whole $\bar{S}$, it stores the minimal portion of information which is able to cooperatively define $\bar{S}$. 
\newline
It is worth noticing that the property of uniquely reconstructing a set via its projections onto a finite number of selected subspaces is not general and holds, as stated in Theorem \ref{thm:bar_S_centralized_from_projection}, for extrusion generated sets associated to a given collection of axis sets. 
\newline
For the networked dynamical system case, the findings of Theorem \ref{thm:bar_S_centralized_from_projection} will be conveniently exploited to derive axis sets, local to each dynamical agents, able to reconstruct the backward reachable set of the overall system. 
Such reconstruction will be performed in a distributed way and not only through the centralized computation in \eqref{eq:bar_S_centralized_from_projection}. 
Still, according to the findings of Theorem \ref{thm:bar_S_centralized_from_projection} the   decomposition of the centralized problem into a distributed one will guarantee no loss of information. 
\end{remark}

The concepts introduced here will be further extended in Section \ref{sec:distributed_extrusion_generated_set}, where we will introduce distributed operations. Therefore, Section \ref{sec:distributed_extrusion_generated_set} can be seen as a continuation of the mathematical background presented here. Nevertheless, we devote for them a new section so as to timely provide, in what follows, the problem statement of this work.

\section{Preliminaries and problem statement}\label{sec:problem_statement}
In this section we firstly provide definitions and concepts 
to introduce the problem formulation and then we give the problem statement.

\subsection{Constrained dynamical system reachability}\label{subsec:dynamical_system_reachability}

Let us first introduce the constrained, discrete-time, dynamical system
\begin{subequations}\label{eq:sys_overall}
\begin{align}
 &x(t+1)=\mathcal{X}(x(t),u(t)),\label{eq:dynamics_overall}\\
 &(x^T(t),u^T(t))^T\in I\subseteq \Bbb{R}^{n+m},\label{eq:mutual_constraints_overall}\\
 & x(t)\in X \subseteq \Bbb{R}^n,\,u(t)\in U\subseteq \Bbb{R}^m. \label{eq:state_input_constraints_overall}
\end{align}
\end{subequations}
where $x$ and $u$ are the state and input vectors, respectively,  and $\mathcal{X}: \Bbb{R}^{n}\times \Bbb{R}^m \to \Bbb{R}^n$.

\begin{remark}\label{rem:explaination_overall_system}
The aim of the paper is to provide a theory for distributed reachability in networked dynamical systems. In view of this, the above system \eqref{eq:sys_overall} is described through a generic nonlinear function $\mathcal{X}(\cdot)$ and possibly static constraints on the states and inputs \eqref{eq:mutual_constraints_overall}.  Both the dynamics and the static constraints are provided in a generic way and no further specification or hypothesis are required. Furthermore, with  \eqref{eq:state_input_constraints_overall} we give constraints on the states and inputs only. Technically, such constraints could be removed and included in \eqref{eq:mutual_constraints_overall}. However, for the sake of clarity in the definition of the problem (as it will be clearer later), we prefer to separate the joint constraints from those related to states and inputs only.  
\end{remark}

We give the following definition.

\begin{defi}\label{def:reachability_H_steps}
Consider the constrained discrete time dynamical system \eqref{eq:dynamics_overall}$\--$\eqref{eq:state_input_constraints_overall}.

Given a set\footnote{
Hereafter, barred sets symbols, e.g., $\bar{S}$, will denote a set obtained with ``centralized'' operations. 
}
$\bar{S}_k\subseteq \bar{P}_k\subseteq X$ and a set $\bar{S}_h\subseteq \bar{P}_h\subseteq X$, we say that $\bar{S}_h$ is {\em reachable} from $\bar{S}_k$ in 
$H\in \Bbb{N}$ steps, iff for any $x(0)\in\bar{S}_k$ there exists a sequence of inputs $u(0), u(1),\dots, u(H-1)\in U$ such that: 
$x(1), \dots, x(H-1)\in \bar{P}_k$; the sequence $(x^T(0),u^T(0))^T, \dots,(x^T(H-1),u^T(H-1))^T \in I$; the terminal state $x(H)\in\bar{S}_h$.

We denote with $\bar{S}_k \rightrsquigarrow^H \bar{S}_h$ the case $\bar{S}_h$ is reachable from~$\bar{S}_k$ in $H$ steps. We denote with $\bar{S}_k \nrightrsquigarrow^H \bar{S}_h$ otherwise. 

For brevity, we will always omit $H$ and we will simply write $\bar{S}_k \rightrsquigarrow \bar{S}_h$ or $\bar{S}_k \nrightrsquigarrow \bar{S}_h$.
\end{defi}

In Definition~\ref{def:reachability_H_steps}, the sets $\bar{P}_k,\bar{P}_h$ are two (among many possible other) sets the state space $X$ is originally partitioned into. The subsets $\bar{S}_k,\bar{S}_h$, as well as the dynamical state trajectory, belong to such sets.
 It is worth mentioning that, in reachability problems where such partition is not of interest, $\bar{P}_k=\bar{P}_h=X$ can be assumed.
On the other hand, such partition turns useful in some specific contexts, such as bisimulation applications (see for example~\cite{nioz:12}), and is therefore included in our reachability problem formulation.

\begin{remark}\label{rem:horizon_H}
Notice that in our analysis we consider a generic $H\geq 0$, with $H\in\Bbb{N}$. Indeed, although for many backward reachability investigations in the literature it is assumed $H=1$, the theory we develop applies to a generic positive integer $H$. In line with the rest of the paper, where we provide general results, we keep the reachability horizon $H$ generic. 
\end{remark}

We also give the following definition. 
\begin{defi}\label{def:pre}
Consider the set $\bar{S}_h\subseteq \bar{P}_h\subseteq X$ and the dynamics \eqref{eq:dynamics_overall}$\--$\eqref{eq:state_input_constraints_overall}. 

The $\mathrm{Pre}(\bar{S}_h)\subseteq \Bbb{R}^n$  with respect to \eqref{eq:dynamics_overall}$\--$\eqref{eq:state_input_constraints_overall} is defined as the set such that $\mathrm{Pre}(\bar{S}_h)\rightrsquigarrow \bar{S}_h$ and $\Bbb{R}^n\backslash\mathrm{Pre}(\bar{S}_h) \nrightrsquigarrow \bar{S}_h$.
\end{defi}

We also give the following definition.

\begin{defi}\label{def:dynamical_flow_overall}
Let us consider system \eqref{eq:sys_overall} and the sets $\bar{S}_k,\,\,\bar{S}_h$, with $\bar{S}_k=\mathrm{Pre}(\bar{S}_h)$. 
Also let us consider vectors  $\mathbf{u}=[u^T(0), \dots, u^T(H)]^T\in U^{H+1}\subseteq \Bbb{R}^{(H+1)m}$, with $U^{H+1}$ the set obtained from the Cartesian product of $U$ with itself $H+1$~times. 
\newline
We define $\bar{\Phi}_{kh}\subseteq \Bbb{R}^{n+(H+1)m}$ the {\em admissible control sequence} of system~\eqref{eq:sys_overall}, that is the set of all the points 
$(x^T(0),\mathbf{u}^T)^T\in \bar{S}_k\times U^{H+1}$ such that, starting from $x(0)$ and applying the sequence of inputs $\mathbf{u}$, point $x(H)$ obtained through the dynamics \eqref{eq:dynamics_overall} is such that $x(H)\in\bar{S}_h$ and constraints \eqref{eq:mutual_constraints_overall}$\--$\eqref{eq:state_input_constraints_overall} are satisfied for $t=0, \dots, H-1$.
\end{defi}
\begin{remark}\label{rem:dynamical_flow_overall}
Notice that, in the above definition, $u(H)$~does not concur in determining $x(H)$ and, therefore, can assume any value. The reason why such input is included in the admissible control sequence definition is only due to a symmetry in the notation, as again it will be clearer later. 
\end{remark}
The problem of calculating the set~$\mathrm{Pre}$ (and related admissible control sequence $\bar{\Phi}_{kh}$) is called {\em backward reachability}, hereafter simply referred to as {\em reachability problem}.

\subsection{Networked local reachability}\label{subsec:dynamical_system_reachability_local}

The aim of this paper is to provide an algorithm for the solution of the reachability problem for {\em large scale networked systems}, that is, for systems with large~$n,~m$ in~\eqref{eq:sys_overall}, characterized by a networked structure as described in what follows.
\newline
Let us consider an ensemble of $N\geq 2$ dynamical systems (also termed {\em agents}), with $N$ a positive finite integer (possibly large). Also, let us suppose without loss of generality that the agents are labelled according to an unique index $i=1,\dots, N$. Such an index will implicitly define an order among them.

The constrained dynamics of each system $i=1,\dots, N$ are expressed by
\begin{subequations}\label{eq:local_dynamics}
\begin{align}
& x_i(t+1)=\mathcal{X}_i(x_i(t),\{x_j(t)\}_{j\in \mathcal{N}_{\mathcal{X},i}},u_i(t),\{u_j(t)\}_{j\in \mathcal{N}_{\mathcal{X},i}}),\label{eq:dynamics_i}\\
&\mathcal{L}_{il}(x_i(t),\{x_j(t)\}_{j\in \mathcal{N}_{I,i}},u_i(t),\{u_j(t)\}_{j\in \mathcal{N}_{I,i}})\sim_{il} 0, \nonumber \\
&l=1,\dots,q_i,\label{eq:mutual_constraints_i}\\
& x_i(t)\in X_i \subseteq \Bbb{R}^{n_i}, \, u_i(t)\in U_i\subseteq\Bbb{R}^{m_i}.  \label{eq:state_input_constraints_i}
\end{align}
\end{subequations}
where $\sim_{il}\in\{<, =, \leq \}$ and $q_i$ is the number of constraints associated with $i$.  The sets $\mathcal{N}_{\mathcal{X},i}$ and $\mathcal{N}_{I,i}$ collects the indices of the other agents 
of the ensemble that are linked to system $i$, respectively through the dynamics and through constraints. 
\begin{defi}\label{def:neighbourhood_i}
We call {\em neighbourhood} of system $i$ the ordered set $\mathcal{N}_i=\mathcal{N}_{\mathcal{X},i}\cup\mathcal{N}_{I,i}\cup \{i\}$, i.e., 
the set of all the agents that interact with $i$ and ordered according to their index. 
Notice that $i$ is a neighbour of itself since it is linked with itself via constraints \eqref{eq:state_input_constraints_i} that, as anticipated in Remark \ref{rem:explaination_overall_system} for the overall system, we conveniently separate from the joint constraints so as to highlight the fact that $i\in\mathcal{N}_i$.
\end{defi}

Now, set $n=\sum_{i=1}^N n_i$ and $m=\sum_{i=1}^N m_i$ and let $x=(x_1^T,\dots,x_N^T)^T\in\Bbb{R}^n$ and $u=(u_1^T,\dots,u_N^T)^T\in\Bbb{R}^m$ be the overall system's state and  input vectors, respectively. 
Let us also define the set 
\begin{align*}
I:=&\left\{(x^T,u^T)^T\in\Bbb{R}^{n+m}: \mathcal{L}_{il}(x,u)\sim_{il} 0, \right.\\
& \left. \forall l=1, \dots, q_i,\,\forall i=1,\dots, N \right\}.
\end{align*}

It is easy to notice that the overall system dynamics can be obtained by stacking the dynamics \eqref{eq:dynamics_i} and the state vectors as $\mathcal{X}(\cdot)=[\mathcal{X}^T_1(\cdot), \dots,\mathcal{X}^T_N(\cdot)]^T$ and setting $X=X_1\times\dots\times X_N$, $U=U_1\times\dots\times U_N$. 
In such a way, the overall system dynamics 
turns into \eqref{eq:sys_overall}.  
For convenience, we refer to the reachability problem for the system \eqref{eq:sys_overall} 
as {\em centralized}. 
In contrast, we will also provide in what follows {\em local} concepts of reachability that will turn to be useful for providing the problem statement of this paper. 
To start with, we consider here vectors $\mathbf{x}_i(t)\in \Bbb{R}^{\sum_{j\in\mathcal{N}_i}n_j}$ and $\mathbf{u}_i(t)\in \Bbb{R}^{\sum_{j\in\mathcal{N}_i}m_j}$. These two vectors respectively stack the states and the inputs of all the agents in $\mathcal{N}_i$ (and so including $i$ itself) at time $t$, ordered according to the same order of $\mathcal{N}_i$. 
\begin{defi}\label{def:local_reachability}
Consider a networked ensemble given by~\eqref{eq:local_dynamics} with $i=1,\dots, N$. 
\newline
Let us focus on an agent $i$ and let us consider a set $S_{k,i}\subseteq P_{k,i}\subseteq \prod_{j\in\mathcal{N}_i}  X_j$ and a set  $S_{h,i}\subseteq P_{h,i}\subseteq \prod_{j\in\mathcal{N}_i}  X_j$. We say that $S_{h,i}$ is {\em locally reachable} from $S_{k,i}$ in $H$~steps iff for any $\mathbf{x}_i(0)\in S_{k,i}$ there exists a sequence of states $\mathbf{x}_i(1),\mathbf{x}_i(2),\dots, \mathbf{x}_i(H)$ and inputs $\mathbf{u}_i(0),\mathbf{u}_i(1),\dots, \mathbf{u}_i(H-1)$ such that $\mathbf{x}_i(H)\in S_{h,i}$ and  the dynamic constraint \eqref{eq:dynamics_i} and the static constraints \eqref{eq:mutual_constraints_i}$\--$\eqref{eq:state_input_constraints_i} are satisfied, for $t=0, \dots, H-1$.
\newline
We denote with $S_{k,i}\rightrsquigarrow_i S_{h,i}$ the case $S_{h,i}$ is locally reachable from $S_{k,i}$. We denote with $S_{k,i}\nrightrsquigarrow_i S_{h,i}$ otherwise.
\end{defi}
Although analogous to Definition \ref{def:reachability_H_steps}, Definition \ref{def:local_reachability} contains some peculiarities we highlight in the following remark.
\begin{remark}\label{rem:local_reachability}
Notice that the definition of local reachability for a single agent $i$ provided in Definition \ref{def:local_reachability} not only involves states $x_i(t)$ and inputs $u_i(t)$ of agent $i$, but also those related to its neighbours in $\mathcal{N}_i$, whose stacks are respectively denoted by $\mathbf{x}_i(t)$ and $\mathbf{u}_i(t)$. Specifically, since the neighbours $\mathcal{N}_i$ directly influence the agent $i$, the local reachability problem involves, in general, both in the starting state region $S_{k,i}$ and the ending state region $S_{h,i}$, all the states in $\mathcal{N}_i$. Similarly, the trajectory of states and input which allows agent $i$ to transit from these two regions must fulfil constraints \eqref{eq:dynamics_i}$\--$\eqref{eq:state_input_constraints_i} which involve all the neighbourhood $\mathcal{N}_i$. 
\newline
Notice also that, from the local perspective of agent $i$ only, states and inputs of all the neighbourhood are seen, in Definition \ref{def:local_reachability}, as decision variables.
\end{remark}
In what follows we provide a definition which represents the local counterpart of Definition \ref{def:pre}.
\begin{defi}\label{def:local_pre}
Consider the networked ensemble~\eqref{eq:local_dynamics}
with $i=1,\dots, N$. 
For agent $i$, consider the set $S_{h,i}\subseteq P_{h,i}\subseteq \prod_{j\in\mathcal{N}_i}  X_j$.  
\newline
The $\mathrm{Pre}_i(S_{h,i})\subseteq \Bbb{R}^{\sum_{j\in\mathcal{N}_i}}$  with respect to \eqref{eq:dynamics_i}$\--$\eqref{eq:state_input_constraints_i} is defined as the set such that $\mathrm{Pre}_i(S_{h,i})\rightrsquigarrow_i S_{h,i}$ and $\Bbb{R}^n \backslash  \mathrm{Pre}_i(S_{h,i}) \nrightrsquigarrow_i S_{h,i}$.
\end{defi}

Having defined $\mathcal{N}_i$ in Definition \ref{def:neighbourhood_i}, the {\em influencing graph}  (for details on graph theory we refer the reader to \cite{goro:01}) of the overall system is defined as $\mathcal{G}=(R,E_{\mathrm{in}})$, with $R=\{1,\dots, N\}$ and $E_{\mathrm{in}}=\{(i,j): j\in\mathcal{N}_i, \,\, i,j\in R\}$. Such a graph is in general directed since, while node $i$ may be subject to the ``influence" of node $j$ according to Definition \ref{def:neighbourhood_i}, the converse may  not be true. Notice also that such a graph may not be connected.

In our paper, we consider the setup where each agent is equipped with local computational and communication capabilities. In particular, we suppose that each system communicates with a subset of other agents according to a {\em communication graph} $\mathcal{G}_c=(R, E)$, where $E=\{(i,j): (i,j)\in E_{\mathrm{in}} \,\, \mathrm{or}\,\, (j,i)\in E_{\mathrm{in}}\}$.
Hence, the communication graph is undirected and built upon the directed graph~$\mathcal{G}$, where a communication channel is established between two nodes, say $i$ and $j$, if either $i$ has a direct influence on $j$ or $j$ has a direct influence on $i$. 

The neighborhood of node $i$ with respect to the communication graph $\mathcal{G}_c$ is denoted by $\mathcal{M}_i$. Note that $i\in \mathcal{M}_i$ and, therefore, $i$ communicates with itself.
Also, we define $\mathcal{M}_i$ as an ordered sets analogously to  $\mathcal{N}_i$ (so inheriting the same order of the agents implicitly provided by their index).

\subsection{Axis sets for the networked reachability problem}\label{subsec:axis_sets_distributed_reachability}
Concepts provided in Section \ref{sec:extrusion_generated_set} are here exploited to define axis sets for the distributed reachability problem. Such axis sets will be both useful in the definition of the problem statement (Section \ref{subsec:problem_statement}) and in the proposed solution for the distributed reachability (Section \ref{sec:distributed_reachability}).

First of all, let us exploit the definitions of $x,u$ given in Section \ref{subsec:dynamical_system_reachability_local}. From these, let us consider the stack\footnote{Variable $u(H)$ will not take part at any of the following reasoning. It is included in the variables stack  only for notational convenience.}
\newline
$\mathbf{z}=[x^T(0), u^T(0), x^T(1), u^T(1),\dots, x^T(H), u^T(H)]^T\in \Bbb{R}^{(H+1)(n+m)}$ , where $m,n$ are given in Section  \ref{sec:problem_statement}.

The stack $\mathbf{z}$ contains, for every time step from $0$ to $H$, the states and inputs of all the networked subsystems.
We define the two sets, depending on the variables $t$ and~$i$, 
\begin{align*}
\tilde{\mathcal{B}}_{x,t,i}:=&\left\{
t(n+m)+ \sum_{j=1}^{i-1}n_j+1, \right. \\
&\left. t(n+m)+\sum_{j=1}^{i-1}n_j+2,
\dots,
t(n+m)+\sum_{j=1}^{i-1}n_j+n_i \right\},\\
\tilde{\mathcal{B}}_{u,t,i}:=&\left\{
t(n+m)+n+ \sum_{j=1}^{i-1}m_j+1, \right. \\
& t(n+m)+n+\sum_{j=1}^{i-1}m_j+2,
\dots, \\
& \left. t(n+m)+n+\sum_{j=1}^{i-1}m_j+m_i \right\},
\end{align*}
where $n_i, m_i$ have been introduced in Section \ref{subsec:dynamical_system_reachability_local}.

The axis sets $\tilde{\mathcal{B}}_{x,t,i}$ and  $\tilde{\mathcal{B}}_{u,t,i}$ are, respectively, the coordinates where the states and the inputs of system $i$ at time $t$ are located in the vector $\mathbf{z}$.

Also, we define $\tilde{\mathcal{B}}_{t,i}:=\tilde{\mathcal{B}}_{x,t,i}\cup\tilde{\mathcal{B}}_{u,t,i}$. 
Thanks to these three sets we define the axis sets related to each subsystem $i$ along with its neighbourhood $\mathcal{M}_i$ at each time instant $t$ as 
\begin{align*}
\mathcal{B}_{x,t,i}:=&\bigcup_{j\in\mathcal{M}_i}\tilde{\mathcal{B}}_{x,t,j},\\
\mathcal{B}_{u,t,i}:=&\bigcup_{j\in\mathcal{M}_i}\tilde{\mathcal{B}}_{u,t,j},\\
\mathcal{B}_{t,i}:=&\mathcal{B}_{x,t,i} \cup \mathcal{B}_{u,t,i},\\
\mathcal{B}_{x,t,\mathcal{M}_i}:=&\bigcup_{j\in\mathcal{M}_i}
\mathcal{B}_{x,t,j},\\
\mathcal{B}_{u,t,\mathcal{M}_i}:=&\bigcup_{j\in\mathcal{M}_i}
\mathcal{B}_{u,t,j},\\
\mathcal{B}_{t,\mathcal{M}_i}:=&\bigcup_{j\in\mathcal{M}_i}
\mathcal{B}_{t,j}.
\end{align*}

Also, we define the axis sets for each node  $i$ for $H$ time steps as
\begin{align*}
\mathcal{B}^H_{x,i}:=& \bigcup_{t=0}^H \mathcal{B}_{x,t,i},
\\
\mathcal{B}^H_{u,i}:=& \bigcup_{t=0}^H \mathcal{B}_{u,t,i},
\\
\mathcal{B}^H_{i}:=& \bigcup_{t=0}^H \mathcal{B}_{t,i},
\\
\mathcal{B}^H_{u,\mathcal{M}_i}:=& \bigcup_{t=0}^H \mathcal{B}_{u,t,\mathcal{M}_i},
\\
\mathcal{B}^H_{\mathcal{M}_i}:=& \bigcup_{t=0}^H \mathcal{B}_{t,\mathcal{M}_i}.
\end{align*}

Let us also define $\bar{\mathcal{B}}^H:=\bigcup_{i=1}^N \mathcal{B}_i^H$, $\bar{\mathcal{B}}_{x,t}:=\bigcup_{i=1}^N \mathcal{B}_{x,t,i}$ and $\bar{\mathcal{B}}_u^H:=\bigcup_{j=1}^N \mathcal{B}_{u,j}^H$.

The meaning of the above defined axis sets may seem not immediate at a first glance. A rough interpretation is provided in the following remark.
\begin{remark}\label{rem:interpretation_axis_sets}
The stack vector  $\mathbf{z}$ contains pairs of the overall system states and inputs (ordered according to their index) for all the time instants $t=0,\dots, H$. The entry indices at which some of such states and inputs are located in vector $\mathbf{z}$ is provided by the axis sets defined before. 
\newline
In particular, set $\tilde{\mathcal{B}}_{x,t,i}$ provides the location indices in $\mathbf{z}$ of $x_i(t)$, that is the state of agent $i$ at time $t$. Analogously, $\tilde{\mathcal{B}}_{u,t,i}$ provides the location of $u_i(t)$. 
Sets $\mathcal{B}_{x,t,i}$ and $\mathcal{B}_{u,t,i}$ respectively provide the location of all the states and inputs of all the neighbours of system $i$ at time $t$, i.e., $x_j(t)$ and $u_j(t)$, with $j\in\mathcal{M}_i$. Notice that such elements are located, in general, not in contiguous positions in  $\mathbf{z}$.  All the other axis sets are suitable unions of such two sets. 
For example, $\mathcal{B}_{t,i}$ represents the indices of the entries in~$\mathbf{z}$ of all the variables affecting system $i$ at time $t$, while  $\mathcal{B}_{t,\mathcal{M}_i}$ represents the indices of all the variables affecting system $i$ neighbourhood at time~$t$. 
\newline
Also, $\mathcal{B}^H_i$ contains the indices of all the variables (states and inputs) affecting agent $i$ whole trajectory, while $\mathcal{B}^H_{\mathcal{M}_i}$ contains the indices of the variables affecting the trajectory of agent $i$ neighbourhood. 
\newline
The set $\bar{\mathcal{B}}^H$ contains all the indices of all the variables of the overall system (i.e., it contains integer numbers from $1$ up to the dimension of $\mathbf{z}$), while $\bar{\mathcal{B}}_{x,t}$ and $\bar{\mathcal{B}}_u^H$ contain, respectively, the location of all the states in the networked ensemble at a specific time instant $t$ and the location of all the inputs of the agents in the ensemble for the whole $H$-steps trajectory. 
\end{remark}

\subsection{Problem statement}\label{subsec:problem_statement}
In this paper, we aim at solving the backward reachability problem for a networked system of equations \eqref{eq:local_dynamics}. To do so, we propose a theory that will allow to solve such a problem in a completely distributed way. 
\newline
First of all, as highlighted in Section \ref{subsec:dynamical_system_reachability_local}, the networked system can be obviously seen as a unique overall dynamical system where a centralized reachability problem can be set so as to compute $\bar{S}_k=\mathrm{Pre}(\bar{S}_h)$ and $\bar{\Phi}_{kh}$ as in Definition \ref{def:pre} and Definition \ref{def:dynamical_flow_overall}, respectively. Notice that, according to the axis sets defined in Section \ref{subsec:axis_sets_distributed_reachability}, we have that $\bar{S}_k\subseteq \Bbb{R}^{|\bar{\mathcal{B}}_{x,0}|}$ and $\bar{\Phi}_{kh}\subseteq \Bbb{R}^{|\bar{\mathcal{B}}_{x,0}\cup \bar{\mathcal{B}}^H_u|}$.
\newline
Another important aspect to notice is that, in the networked setting, the terminal set $\bar{S}_h$ can be seen as an extrusion generated set obtained according to 
\begin{equation}\label{eq:bar_S_h_from_S_h_i}
\bar{S}_h=\bigcap_{i=i}^N\mathcal{E}_{\mathcal{B}_{x,H,i}}^{\bar{\mathcal{B}}_{x,H}}(S_{h,i}).
\end{equation}
\newline
Indeed, as also highlighted in Remark \ref{rem:local_reachability},
for each agent~$i$ joint constraints only exist among those agents in~$\mathcal{N}_i$ (and so also in $\mathcal{M}_i$). This includes the terminal set region and, therefore, the overall terminal state $x(H)$ must satisfy the relation $\mathcal{P}_{\mathcal{B}_{x,H,i}}^{\bar{\mathcal{B}}_{x,H}}(x(H))\in S_{i,h}$ for all $i=1,\dots, N$ and, therefore, $x(H)\in\bar{S}_h$ with $\bar{S}_h$ defined in \eqref{eq:bar_S_h_from_S_h_i}.
\newline
In this paper, we will show that a relation analogous to~\eqref{eq:bar_S_h_from_S_h_i} is satisfied also by set $\bar{S}_k$ which will be shown being an extrusion generated set associated to the collection~$\mathcal{B}_{x,0,i}$ with $i=1,\dots, N$. Specifically, taking into account Theorem \ref{thm:bar_S_centralized_from_projection}, $\bar{S}_k$ can be generated by the minimal set of information distributed at each node level $i$ as sets $\bar{S}_{k,i}$, with $i=1,\dots, N$ and 
\begin{equation}\label{eq:bar_S_k_i_from_bar_S_k}
\bar{S}_{k,i}=\mathcal{P}_{\mathcal{B}_{x,0,i}}^{\bar{\mathcal{B}}_{x,0}}(\bar{S}_k).
\end{equation}
\newline
Such sets $\bar{S}_{k,i}$ can be obviously found via first solving in a centralized way the reachability problem of the overall system as a whole in order to determine $\bar{S}_k$. Then, the $\bar{S}_{k,i}$ can be obtained through \eqref{eq:bar_S_k_i_from_bar_S_k}. Both the operations are, however, centralized and involve all the decision variables of the overall system\footnote{This aspect holds irrespectively of the algorithm that might be used for computing the reachability.}, making them impractical for not trivial size network dimension.
\newline
The objective of the paper is to develop a theoretical framework able to shield light on the intrinsic link between the concept of backward reachability for large scale systems and its distributed counterpart, providing a systematic way to convert a general nonlinear constrained networked problem into a distributed one. Specifically, the paper will propose a way to determine, at the level of each agent $i$, the set $\bar{S}_{k,i}$ via only solving a local backward reachability problem as in Definition \ref{def:local_reachability}, starting from $\bar{S}_{h,i}$ and employing local information exchanges with the agent $i$ neighbours. From a networked viewpoint, $\bar{S}_k$ will implicitly be coded at the network level through the $\bar{S}_{k,i}$'s stored at node level without any centralized operation. 
\newline
In addition to this, for the backward reachability problem of the overall system,  the admissible control sequence $\bar{\Phi}_{kh}$ as defined in Definition \ref{def:dynamical_flow_overall} might also be of interest since it contains the control sequences that allow, in $H$ steps, to move the network from the global state space region $\bar{S}_k$ to the target region $\bar{S}_h$. In this paper we will show that set $\bar{\Phi}_{kh}$ is  also an extrusion generated set with respect to the (local) collection of axis sets $\mathcal{B}_{x,0,i}\cup \mathcal{B}_{u,i}^H$ and that can therefore be reconstructed through the local information $\bar{\Phi}_{kh,i}$'s at node level, where 
\begin{equation}\label{eq:bar_Phi_kh_i_from_bar_Phi_kh}
\bar{\Phi}_{kh,i}=\mathcal{P}_{\mathcal{B}_{x,0,i}\cup \mathcal{B}_{u,i}^H}^{\bar{\mathcal{B}}_{x,0}\cup \bar{\mathcal{B}}_{u}^H}(\bar{\Phi}_{kh}).
\end{equation}
Analogously to $S_{k,i}$, also for $\bar{\Phi}_{kh,i}$ we will provide a completely decentralized way for its computation without passing through the centralized \eqref{eq:bar_Phi_kh_i_from_bar_Phi_kh}.
\newline
From a broader viewpoint, for a nonlinear constrained system \eqref{eq:sys_overall} where state and input variables form a networked pattern so as to be rewritten according to \eqref{eq:local_dynamics}, the paper in hand shows how the backward reachability problem can be rewritten in an equivalent (i.e., with no information losses) distributed way.
The theory here proposed is completely abstract and applies to any kind of system dynamics, constraints and topology of the state/input space irrespectively on how an analytical or numerical solution can be (if possible) actually computed and irrespectively of the solver that might be used.

\section{Mathematical framework for distributed extrusion generated set}\label{sec:distributed_extrusion_generated_set}
In this section, we propose some important results which will ultimately allow to solve the problem described in Section \ref{subsec:problem_statement}. To do so, we put again the focus on  the mathematical concepts given in Section \ref{sec:extrusion_generated_set} and further extend them. Specifically, we provide a distributed algorithm able to compute, via iterative local operations and local information exchanges, the projections of an extrusion generated set (which will never need to be computed) defined as  in   \eqref{eq:projection_extrusion_generated_set}.
\newline 
It is worth noticing that this section can be seen as a continuation of the mathematical background part of Section \ref{sec:extrusion_generated_set} (interrupted so as to timely provide the problem statement of the paper). All the results are developed here in general terms and are instrumental to the derivation of the main results of the paper in Section~\ref{sec:distributed_reachability}.
As such, mathematical objects like the axis sets $\mathcal{B}_i$ considered here are, therefore, not directly related to the distributed reachability problem nor they are
linked here to dynamical systems. Such a link will be provided in Section \ref{sec:distributed_reachability}.

To start with, it is worth highlighting that computing the set \eqref{eq:extrusion_generated_set} requires the extrusion of each set $S_i$ towards all the other axis sets and then a centralized intersection operation, which turns out to be cumbersome when large $N$ and the nontrivial case of identical $\mathcal{B}_i$ is considered. 
In what follows we propose a distributed iterative computation of the extrusion generated set $\bar{S}$. This will turn useful in Section \ref{sec:distributed_reachability}, where the distributed reachability problem will be solved.

Consider $N$ non empty axis sets $\mathcal{B}_1,\mathcal{B}_2,\dots,\mathcal{B}_N$, let again define $\bar{\mathcal{B}}$ as in \eqref{eq:B_bar}, and construct an undirected graph~$\mathcal{G}_c=(R,E)$ such that each node $i\in R$ is associated to~$\mathcal{B}_i$ and there exists an undirected arc $(i,j)\in E$ if the two axis sets associated to nodes $i$ and $j$ are not disjoint, i.e.,~$\mathcal{B}_i\cap \mathcal{B}_j \neq \emptyset$.
The neighbourhood of node $i$ according to $\mathcal{G}_c$ is denoted by $\mathcal{M}_i$. Note that, since $\mathcal{B}_i$ is not disjoint with itself, $i \in \mathcal{M}_i$. Note also that $\mathcal{G}_c$ may not necessarily be a connected graph.

Example \ref{exa:graph} shows the construction of sets $\mathcal{M}_i$ (and, therefore, of graph $\mathcal{G}_c$) from sets $\mathcal{B}_i$.
\begin{remark}
Notice that the graph defined here among axis sets is different from that on of Section \ref{subsec:dynamical_system_reachability_local} involving dynamical systems. 
\newline
The reason why we adopt the same symbol $\mathcal{M}_i$ is intentional and strongly linked with the axis sets introduced in Section \ref{subsec:dynamical_system_reachability_local}.
Indeed, via such axis sets, we will be able to link the mathematical concepts we will develop in this section to the dynamical system case. For example, considering a graph where at each node $i$ we locate the corresponding $\mathcal{B}_i^H$, it is immediate to notice that the neighbourhood $\mathcal{M}_i$ defined as before, i.e.  the set of $j$ such as $\mathcal{B}_i^H\cap \mathcal{B}_j^H\neq \emptyset$, provides the same set as the one in Section~\ref{subsec:dynamical_system_reachability_local}. 
\newline
In general, for dynamical systems, we will always resort in the rest of the paper to the axis set formalism and related operations. 
\end{remark}

As initialization phase of the algorithm, each node exchanges its $\mathcal{B}_i$ with its neighbors. Such information exchange is necessary to execute Algorithm \ref{alg:distributed_extrusion_generated_set}, which runs synchronously for all the nodes at the time instants $\kappa$. Such time instants are not related to the time variable~$t$ in \eqref{eq:dynamics_i}$\--$\eqref{eq:state_input_constraints_i}, but rather counts the clock tick when the local processing occurs at all nodes.
The Algorithm \ref{alg:distributed_extrusion_generated_set} calculates in a distributed way the set $\bar{S}=\bigcap_{i=1}^N \mathcal{E}_{\mathcal{B}_i}^{\bar{\mathcal{B}}}(S_{i,0})$ and its projections $\bar{S}_i$ (as defined in \eqref{eq:projection_extrusion_generated_set}), where, for each node $i$, the sets $S_i(0)=S_{i,0}\subseteq \Bbb{R}^{|\mathcal{B}_i|}$ initialize the algorithm.

\begin{algorithm}
\caption{Distributed extrusion generated set: node~$i$ projection}\label{alg:distributed_extrusion_generated_set}

\begin{algorithmic}[1]

\ForAll{$\kappa\geq 0$}

\State Send $S_i(\kappa)$ to all $j\in\mathcal{M}_i$;

\State Receive $S_j(\kappa)$ from all $j\in\mathcal{M}_i$;

\State Compute 
\begin{equation}\label{eq:distributed_extrusion_generated_set}
S_i(\kappa+1)=\mathcal{P}_{\mathcal{B}_i}^{\mathcal{B}_{\mathcal{M}_i}}\left[
\bigcap_{j\in\mathcal{M}_i}
\mathcal{E}_{\mathcal{B}_j}^{\mathcal{B}_{\mathcal{M}_i}}\left(S_j(\kappa)\right)
\right],
\end{equation}

with $\mathcal{B}_{\mathcal{M}_i}=\bigcup_{j\in\mathcal{M}_i}\mathcal{B}_j$;

\State Set $\kappa\leftarrow \kappa+1$;

\EndFor

\end{algorithmic}
\end{algorithm}

Next in this section, we will show that successive iterations of \eqref{eq:distributed_extrusion_generated_set} refine the estimates~$S_i(\kappa+1)$ of $\bar{S}_i$ 
thus ``converging'' to $\bar{S}_i$. 
This will be done through a series of intermediate results. 

\begin{lemma}\label{lem:succession_nested_inclusions}
Consider the $N$ axis sets $\mathcal{B}_1,\mathcal{B}_2,\dots,\mathcal{B}_N$ and the sets $S_{i,0}\subseteq \Bbb{R}^{|\mathcal{B}_i|}$, for $i=1,\dots,N$. Also, let $\bar{S}$ be the extrusion generated set, obtained by applying \eqref{eq:extrusion_generated_set} to the sets~$S_{i,0}$, and $\bar{S}_i$ be its projections according to \eqref{eq:projection_extrusion_generated_set}. From Algorithm \ref{alg:distributed_extrusion_generated_set} the two following results hold:
\begin{itemize}
\item[i.] $\bigcap_{i=1}^N\mathcal{E}_{\mathcal{B}_i}^{\bar{\mathcal{B}}}(S_i(\kappa))=\bar{S},\,\,\forall \kappa\geq 0$;
\item[ii.] $\bar{S}_i\subseteq S_i(\kappa)\subseteq S_i(\kappa-1)\subseteq \dots \subseteq S_i(1)\subseteq S_i(0)$ for any $\kappa>1$.
\end{itemize}
\end{lemma}

Let us now define, for all $i=1,\dots, N$,
\begin{equation}\label{eq:S_i_star}
S_i^{*}=\bigcap_{\kappa=0}^{+\infty} S_i(\kappa).
\end{equation}

The following lemma holds. 

\begin{lemma}\label{lem:fixed_point}
Let us consider Algorithm \ref{alg:distributed_extrusion_generated_set} and the sets $S_i^{*}$ as defined in \eqref{eq:S_i_star}. The following results hold:
\begin{itemize}
\item [i.] $S_i^{*}= \mathcal{P}_{\mathcal{B}_i}^{\mathcal{B}_{\mathcal{M}_i}}\left[\bigcap_{j\in \mathcal{M}_i}\mathcal{E}_{\mathcal{B}_j}^{\mathcal{B}_{\mathcal{M}_i}} (S_j^{*})\right], \,\, \forall i$;

\item[ii.] $S^{*}=\bar{S}$, with $S^{*}$ defined as $S^{*}:=\bigcap_{i=1}^{N}\mathcal{E}_{\mathcal{B}_i}^{\bar{\mathcal{B}}}(S_i^{*})$.
\end{itemize}
\end{lemma}

Lemma \ref{lem:fixed_point} states the existence of a ``fixed point'' for Algorithm \ref{alg:distributed_extrusion_generated_set}. However, it remains to prove that such fixed point corresponds to the ``tightest'' one, i.e, $S_i^*=\bar{S}_i$.
This is provided in the following theorem.

\begin{theorem}\label{thm:thm_distributed_extrusion_algorithm}
Let us consider Algorithm \ref{alg:distributed_extrusion_generated_set} and the sets $\bar{S}_i$ defined in  \eqref{eq:projection_extrusion_generated_set}, for $i=1,\dots, N$. Then, for all $i=1,\dots, N$ the sequence of sets $\{S_i(\kappa)\}_{\kappa=0}^{+\infty}$ {\em approaches} $\bar{S}_i$ for $\kappa\rightarrow +\infty$. That is, $\bar{S}_i\subseteq \dots \subseteq S_i(\kappa)\subseteq S_i(\kappa-1)\subseteq \dots \subseteq S_i(0)$ with $\bar{S}_i$ fixed point of Algorithm \ref{alg:distributed_extrusion_generated_set}. 
\end{theorem}

\begin{remark}
Via Theorem \ref{thm:thm_distributed_extrusion_algorithm} we show that it is possible to compute, at each node $i$, the projection $\bar{S}_i$ of the centralized extrusion generated set $\bar{S}$ only via iterative local computations. Such a set, which is equal to $S_i^*$ defined in \eqref{eq:S_i_star} is obtained from Algorithm \ref{alg:distributed_extrusion_generated_set} for $\kappa \rightarrow +\infty$.  In what follows, we provide examples where $\bar{S}_i$ is obtained in a finite number of steps. However, in general, it not possible to say that there exists a finite $\bar{\kappa}_i$ such that $\bar{S}_i(\kappa^\prime)=\bar{S}_i(\bar{\kappa}_i)$ for any $\kappa^\prime\geq\bar{\kappa}_i$.
In other word set $\bar{S}_i$ is in general achieved after an infinite number of iterations. 
From a practical viewpoint, this does not constitute a problem since sequence $\bar{S}_i\subseteq \dots \subseteq S_i(\kappa+1)\subseteq S_i(\kappa)\subseteq \dots\subseteq S_i(0)$ refines the solution set $\bar{S}_i$ from outside at any iteration. In numerical setting it is therefore possible to establish termination criteria based on, e.g., a measure of set $S_i(\kappa)$ so as to stop the iterations when a desired precision is reached among successive iterations. 
\newline
Since in this paper we focus on a purely theoretical contribution, we do not provide here further implementation details on such aspects and refer to a future work (a discussion is later provided in Section \ref{sec:discussion}).
\end{remark}

To illustrate Algorithm \ref{alg:distributed_extrusion_generated_set}, we firstly consider the toy Example \ref{exa:extrusion_generated_set} on finite countable sets $S_i$. Also, to further provide insights on how  Algorithm \ref{alg:distributed_extrusion_generated_set} works, a numerical example (Example \ref{exa:extrusion_generated_polyhedron}) is provided, where the projections of a high-dimensional extrusion generated polytope are computed in a distributed way.

\section{Distributed reachability}\label{sec:distributed_reachability}
In this section we exploit Algorithm \ref{alg:distributed_extrusion_generated_set} to solve the distributed  reachability problem for networked systems \eqref{eq:dynamics_i}$\--$\eqref{eq:state_input_constraints_i}.

First, let us consider the sets $\bar{S}_k\subseteq \bar{P}_k$ and $\bar{S}_h\subseteq \bar{P}_h$. In order to study the reachability problem $\bar{S}_k  \rightrsquigarrow \bar{S}_h$ in $H$ steps under the dynamics \eqref{eq:dynamics_i}$\--$\eqref{eq:state_input_constraints_i}, we consider the system in the unknown $\mathbf{z}\in\Bbb{R}^{(H+1)(n+m)}$
\begin{equation}\label{eq:system_of_constraints_overall}
\begin{cases}
\mathcal{P}_{\tilde{\mathcal{B}}_{x,t+1,i}}^{\bar{\mathcal{B}}^H}(\mathbf{z})=\mathcal{X}_i\left(\mathcal{P}_{\mathcal{B}_{t,i}}^{\bar{\mathcal{B}}^H}(\mathbf{z}) \right),\,\, t=0,\dots, H-1, \\
\mathcal{L}_{il}\left(\mathcal{P}_{\mathcal{B}_{t,i}}^{\bar{\mathcal{B}}^H}(\mathbf{z})\right)\sim_{il} 0,\, l=1,\dots,q_i,\,\, t=0,\dots, H-1,\\
\mathcal{P}_{\tilde{\mathcal{B}}_{x,t,i}}^{\bar{\mathcal{B}}^H}(\mathbf{z})\in X_i,\,\,
\mathcal{P}_{\tilde{\mathcal{B}}_{u,t,i}}^{\bar{\mathcal{B}}^H}(\mathbf{z})\in U_i,\,\, t=0,\dots, H,\\
\mathcal{P}_{\bar{\mathcal{B}}_{x,0}}^{\bar{\mathcal{B}}^H}(\mathbf{z})\in \bar{S}_{k},\\
\mathcal{P}_{\bar{\mathcal{B}}_{x,t}}^{\bar{\mathcal{B}}^H}(\mathbf{z})\in \bar{P}_{k},\,\, t=0, \dots, H-1, \\
\mathcal{P}_{\bar{\mathcal{B}}_{x,H}}^{\bar{\mathcal{B}}^H}(\mathbf{z})\in \bar{S}_{h}\subseteq \bar{P}_h,\\
i=1,\dots,N, \\
\end{cases}
\end{equation}
where we recall that $\sim_{il}\in\{<, =, \leq \}$.

Let us call $\bar{S}_{kh}\subseteq \Bbb{R}^{(H+1)(n+m)}$ the {\em solution} of the system \eqref{eq:system_of_constraints_overall} (such solution, in general, does not have an explicit form). 
Roughly speaking, each point of $\bar{S}_{kh}$ represents for all the systems  \eqref{eq:dynamics_i}$\--$\eqref{eq:state_input_constraints_i} in the stack $i=1, \dots, N$ an initial condition and an input admissible trajectory such that the state trajectory obtained from the initial condition under the dynamics \eqref{eq:dynamics_i} satisfies constraints \eqref{eq:mutual_constraints_i} and the state and input constraints satisfy \eqref{eq:state_input_constraints_i}.

The following lemma holds.
\begin{lemma}\label{lem:reachability_overall_system}
Let us consider the system \eqref{eq:system_of_constraints_overall}, and let $\bar{S}_{kh}\subseteq \Bbb{R}^{(H+1)(n+m)}$ be its solution.  Also, let us consider the overall dynamical system \eqref{eq:dynamics_overall}\--\eqref{eq:state_input_constraints_overall} obtained by stacking \eqref{eq:dynamics_i}\--\eqref{eq:state_input_constraints_i}, $i=1,\dots, N$ as described in Section \ref{alg:distributed_Pre_computation}.
Then $\bar{S}_h$ is reachable from $\bar{S}_k$ in $H$ steps if and only if $\mathcal{P}_{\bar{\mathcal{B}}_{x,0}}^{\bar{\mathcal{B}}^H}(\bar{S}_{kh})=\bar{S}_k$.
\end{lemma}

A similar lemma can be given for Definition \ref{def:pre}.
\begin{lemma}\label{lem:pre_overall_system}
Let us consider the system in $\mathbf{z}\in\Bbb{R}^{(H+1)(n+m)}$ obtained from \eqref{eq:system_of_constraints_overall} by removing the constraint $\mathcal{P}_{\bar{\mathcal{B}}_{x,0}}^{\bar{\mathcal{B}}^H}(\mathbf{z})\in \bar{S}_{k}$ and let $\bar{S}_{kh}$ be its solution. 
Also, let us consider the overall dynamical system \eqref{eq:dynamics_overall}\--\eqref{eq:state_input_constraints_overall} obtained by stacking \eqref{eq:dynamics_i}\--\eqref{eq:state_input_constraints_i}, $i=1,\dots, N$ as described in Section \ref{alg:distributed_Pre_computation}.
Then $\mathrm{Pre}(\bar{S}_h)=\mathcal{P}_{\bar{\mathcal{B}}_{x,0}}^{\bar{\mathcal{B}}^H}(\bar{S}_{kh})$.
\end{lemma}

Let us now introduce, for $i=1,\dots, N$, the systems in the unknown $z\in\Bbb{R}^{|\mathcal{B}_i^H|}$
\begin{equation}\label{eq:system_of_constraints_i}
\begin{cases}
\mathcal{P}_{\tilde{\mathcal{B}}_{x,t+1,i}}^{\mathcal{B}_i^H}(z)=\mathcal{X}_i\left(\mathcal{P}_{\mathcal{B}_{t,i}}^{\mathcal{B}_i^H}(z)\right),\,\, t=0,\dots, H-1, \\
\mathcal{L}_{il}\left(\mathcal{P}_{\mathcal{B}_{t,i}}^{\mathcal{B}_i^H}(z)\right)\sim_{il} 0,\, l=1,\dots,q_i,\,\, t=0,\dots, H-1,\\
\mathcal{P}_{\tilde{\mathcal{B}}_{x,t,i}}^{\mathcal{B}_i^H}(z)\in X_i,\,\,
\mathcal{P}_{\tilde{\mathcal{B}}_{u,t,i}}^{\mathcal{B}_i^H}(z)\in U_i,\,\, t=0,\dots, H,\\
\mathcal{P}_{\mathcal{B}_{x,0,i}}^{\mathcal{B}_i^H}(z)\in S_{k,i}\subseteq \Bbb{R}^{|\mathcal{B}_{x,0,i}|},\\
\mathcal{P}_{\mathcal{B}_{x,t,i}}^{\mathcal{B}_i^H}(z)\in P_{k,i}\subseteq \Bbb{R}^{|\mathcal{B}_{x,t,i}|},\,\, t=0, \dots, H-1, \\
\mathcal{P}_{\mathcal{B}_{x,H,i}}^{\mathcal{B}_i^H}(z)\in S_{h,i}\subseteq P_{h,i}\subseteq \Bbb{R}^{|\mathcal{B}_{x,H,i}|},\\
\end{cases}
\end{equation}
with $\bar{S}_k=\cap_{i=1}^N\mathcal{E}_{\mathcal{B}_{x,0,i}}^{\bar{\mathcal{B}}_{x,0}}(S_{k,i})$ and similarly~ 
$\bar{P}_k=\cap_{i=1}^N\mathcal{E}_{\mathcal{B}_{x,,i}}^{\bar{\mathcal{B}}_{x,0}}(P_{k,i})$, 
$\bar{S}_h=\cap_{i=1}^N\mathcal{E}_{\mathcal{B}_{x,H,i}}^{\bar{\mathcal{B}}_{x,H}}(S_{h,i})$, 
\newline 
$\bar{P}_h=\cap_{i=1}^N\mathcal{E}_{\mathcal{B}_{x,H,i}}^{\bar{\mathcal{B}}_{x,H}}(P_{h,i})$.
Note that $|\mathcal{B}_{x,0,i}|=|\mathcal{B}_{x,t,i}|=|\mathcal{B}_{x,H,i}|, \,\,\forall t=0,\dots, H$.

Let us call $S_{kh,i}\subseteq \Bbb{R}^{|\mathcal{B}_i^H|}$  the solution of  \eqref{eq:system_of_constraints_i}. 
Notice that the  structure of $z$ is implicitly defined by the operators $\mathcal{P}_{(\cdot)}^{\mathcal{B}_i^H}$ in the system \eqref{eq:system_of_constraints_i} such that all the entries in $z$ play their role of unknowns within the equations correspondingly to equations in \eqref{eq:dynamics_i}\--\eqref{eq:state_input_constraints_i}.
Notice also that system \eqref{eq:system_of_constraints_i} solves the local reachability relation $S_{k,i}\nrightrsquigarrow_i S_{h,i}$ as per Definition \ref{def:local_reachability}.

We now give the following lemma.

\begin{lemma}\label{lem:bar_S_kh_extrusion_generated}
Let us consider system \eqref{eq:system_of_constraints_overall} and systems~\eqref{eq:system_of_constraints_i}, for $i=1,\dots, N$. Then $\bar{S}_{kh}$ is an extrusion generated set obtained by
\begin{equation}\label{eq:bar_S_kh_extrusion_generated}
\bar{S}_{kh}=\bigcap_{i=1}^N\mathcal{E}_{\mathcal{B}_i^H}^{\bar{\mathcal{B}}^H}(S_{kh,i}).
\end{equation}
\end{lemma}

Formula \eqref{eq:bar_S_kh_extrusion_generated} expresses $\bar{S}_{kh}$ as a centralized extrusion generated set obtained from the $S_{kh,i}$. Therefore, the results obtained in Section \ref{sec:distributed_extrusion_generated_set} can directly be applied. We give the following theorem.

\begin{theorem}\label{thm:distributed_computation_bar_S_kh}
Let us consider system \eqref{eq:system_of_constraints_overall} and systems \eqref{eq:system_of_constraints_i} for all $i=1,\dots, N$. Also, 
let us set $S_{kh,i,0}:=S_{kh, i}$ and let us consider the initialization $S_{kh,i}(0):=S_{kh,i,0}$. The extrusion generated set $\bar{S}_{kh}$ can be computed in a distributed way via
applying  Algorithm \ref{alg:distributed_extrusion_generated_set} particularizing formula \eqref{eq:distributed_extrusion_generated_set}
with the following one
\begin{equation}\label{eq:distributed_extrusion_generated_set_no_disturbances}
S_{kh,i}(\kappa+1)=\mathcal{P}_{\mathcal{B}_i^H}^{\mathcal{B}_{\mathcal{M}_i}^H}\left[
\bigcap_{j\in\mathcal{M}_i}
\mathcal{E}_{\mathcal{B}_j^H}^{\mathcal{B}_{\mathcal{M}_i}^H}\left(S_{kh,j}(\kappa)\right)
\right].
\end{equation}
\end{theorem}

Theorem \ref{thm:distributed_computation_bar_S_kh} establishes a relation among the distributed $S_{kh,i}$ and  $\bar{S}_{kh}$. On the other hand, exploiting  Lemma \ref{lem:reachability_overall_system} and Lemma \ref{lem:pre_overall_system} we have that $\bar{S}_k=\mathrm{Pre}(\bar{S}_h)=\mathcal{P}_{\bar{\mathcal{B}}_{x,0}}^{\bar{\mathcal{B}}}(\bar{S}_{kh})$. 
It is so intriguing to investigate the possibility of computing $\bar{S}_k$ in a distributed way as well as done for $\bar{S}_{kh}$. The following theorem provides a positive answer to such a question.

\begin{theorem}\label{thm:distributed_Pre}
Let us consider the same setup of Theorem \ref{thm:distributed_computation_bar_S_kh}  and let us compute $\bar{S}_{kh,i}$  as fixed point of the distributed Algorithm \ref{alg:distributed_extrusion_generated_set} casted according to \eqref{eq:distributed_extrusion_generated_set_no_disturbances}. Also, let us define
\begin{equation*}
\underline{S}_{k,i}:=\mathcal{P}_{\mathcal{B}_{x,0,i}}^{\mathcal{B}^H_i }(\bar{S}_{kh,i}), 
\end{equation*}
and the set $\bar{S}_{k,i}:=\mathcal{P}_{\mathcal{B}_{x,0,i}}^{\bar{\mathcal{B}}_{x,0}}(\bar{S}_k)$.
Then $\bar{S}_{k,i}=\underline{S}_{k,i}$. 
\end{theorem}

A result similar to Theorem \ref{thm:distributed_Pre} can be considered for the admissible control sequence $\bar{\Phi}_{kh}$ of the overall system as provided in Definition \ref{def:dynamical_flow_overall}. Indeed, it is immediate to notice that $\bar{\Phi}_{kh}$ can also be written as
\begin{equation}\label{eq:flow_overall}
\bar{\Phi}_{kh}:=\mathcal{P}_{\bar{\mathcal{B}}_{x,0}\cup \bar{\mathcal{B}}_{u}^H}^{\bar{\mathcal{B}}^H}(\bar{S}_{kh}),
\end{equation}
with $\bar{\Phi}_{kh,i}$ its projection as per equation \eqref{eq:bar_Phi_kh_i_from_bar_Phi_kh}.  We define
\begin{equation}\label{eq:flow_i}
\underline{\Phi}_{kh,i}:=\mathcal{P}_{\mathcal{B}_{x,0,i}\cup \mathcal{B}_{u,i}^H}^{\mathcal{B}^H_i}(\bar{S}_{kh,i}).
\end{equation}

The following theorem holds. 
\begin{theorem}\label{thm:distributed_flow}
Let us consider the same setup as in Theorem \ref{thm:distributed_computation_bar_S_kh}  and let us compute $\bar{\Phi}_{kh}$, $\bar{\Phi}_{kh,i}$ and $\underline{\Phi}_{kh,i}$ according to \eqref{eq:flow_overall}, \eqref{eq:bar_Phi_kh_i_from_bar_Phi_kh} and \eqref{eq:flow_i}, respectively. Then $\bar{\Phi}_{kh,i}=\underline{\Phi}_{kh,i}$. 
\end{theorem}

The results developed in this section allow to compute in a distributed way the $\bar{S}_k=\mathrm{Pre}(\bar{S}_h)$ under the constrained multi-agent system dynamics \eqref{eq:dynamics_i}$\--$\eqref{eq:state_input_constraints_i} for $\bar{S}_h$ extrusion generated target set. Specifically, each agent computes the admissible states $\bar{S}_{k,i}$ and the associated admissible control inputs $\bar{\Phi}_{kh,i}$. Notice that no loss of information happens when shifting from the centralized to the distributed problem.
The procedure is formalized in Algorithm \ref{alg:distributed_Pre_computation}.

\begin{algorithm}
\caption{Distributed $\bar{S}_{k,i}$ and $\bar{\Phi}_{kh,i}$ computation}\label{alg:distributed_Pre_computation}

\begin{algorithmic}[1]

\State Solve system \eqref{eq:system_of_constraints_i};

\State Compute $S_{kh,i}$ according to \eqref{eq:system_of_constraints_i};

\State Set $S_{kh,i,0}:=S_{kh,i}$;

\State Exchange $S_{kh,i,0}$ with neighbours $j\in\mathcal{M}_i$;



\State Receive $S_{kh,j,0},\,\,\forall j\in\mathcal{M}_i-\{i\}$ from the neighbours, and set initial conditions $S_{kh,j}(0)=S_{kh,j,0}$;

\State Run Algorithm \ref{alg:distributed_extrusion_generated_set} via local exchanges  of $S_{kh,j}(\kappa)$, with $\kappa\geq 0$ and $\forall j\in\mathcal{M}_i$.

The algorithm will provide $\bar{S}_{kh,i}$ as output;

\State Compute  $\bar{S}_{k,i}=\mathcal{P}_{\mathcal{B}_{x,0,i}}^{\mathcal{B}_i^H}(\bar{S}_{kh,i})$;

\State Compute  $\bar{\Phi}_{kh,i}=\mathcal{P}_{\mathcal{B}_{x,0,i}\cup \mathcal{B}_{u,i}^H}^{\mathcal{B}_i^H}(\bar{S}_{kh,i})$.

\end{algorithmic}
\end{algorithm}

\section{Discussion and applicability}\label{sec:discussion}
The aim of this paper is to develop a theory to systematically link centralized vs distributed backward reachability control problem. The manuscript investigates this point from a purely theoretical perspective so as to come up with general conclusions irrespectively of the systems dynamics, types of constraints or state and input region topology. 
\newline
A possible question may be related to how actually solve system \eqref{eq:system_of_constraints_i}. 
Notice that, solving \eqref{eq:system_of_constraints_i} is exactly solving a local reachability problem as per Definition \ref{def:local_reachability}, where for system $i$ the other neighbours' variables are extra decision variables. Such a problem is already the type of problem addressed in several papers dealing with reachability (see the Introduction for a non exhaustive overview on reachability results), where tools exist for particular nonlinear functions and starting/target space sets. The theory developed in this paper is about the fact that, providing that such system can be solved in a numerical/analytical way, networked problem can be solved too via the proposed decomposition whatever size of the problem is considered. 
\newline
Having not restricted the results to any particular nonlinear form, resolution techniques for \eqref{eq:system_of_constraints_i} cannot be specified  here (some relevant approaches are reported in the Introduction). Similarly, the numerical computation of projection and extrusion operators depends on the specific system considered and needs proper consequent numerical investigation, possibly exploiting  also specific peculiarities of the particular class of the investigated systems.  
For example, the relevant case of interconnected linear affine systems (also subjected to disturbance) and with polytopic constraints can be easily addressed through linear programming (LP) operations. Linear affine systems have been highly investigated in the control community due to their adoption in several applications and, furthermore, due to the fact that smooth nonlinear systems may be approximated with piecewise linear affine systems via hybridization/linearization techniques.
\newline
In what follows, we briefly provide an idea on how to cast the theory for this case, while a thorough study will be provided in a dedicated paper. 
First of all, dynamics \eqref{eq:dynamics_i} take the form 
\begin{equation}\label{eq:dynamics_i_linear}
x_i(t+1)=\sum_{j=1}^N A_{ij}x_j(t)+\sum_{j=1}^N B_{ij}u_j(t)+K_i+E_id_i(t),  
\end{equation}
with matrices $A_{ij}$, $B_{ij}$ and $E_i$ of suitable dimensions, and with disturbance $d_i(t)\in D_i\subset \Bbb{R}^{v_i}$ belonging to a bounded polytope. Also, $X_i$ and $U_i$ in \eqref{eq:state_input_constraints_i} are polytopes while for \eqref{eq:mutual_constraints_i} inequality constraints $\mathcal{L}_{il}(\cdot)\leq 0$ will be adopted, with $\mathcal{L}_{il}(\cdot)$ linear affine expressions.  
Notice that, in case $j\notin \mathcal{M}_i$, then $A_{ij}=O_{n_i\times n_j}$ and $B_{ij}=O_{n_i\times m_j}$, with $O$ being the matrix with all null entries.
\newline
Finally, to conclude the setup for the linear affine case, the study of the reachability coded in system 
\eqref{eq:system_of_constraints_overall}, with the help of local systems \eqref{eq:system_of_constraints_i}, will be carried out for the case of polytopes $S_{k,i},P_{k,i},S_{h,i},P_{h,i}$, which obviously lead to extrusion generated polytope sets $\bar{S}_k,\bar{P}_k,\bar{S}_h,\bar{P}_h$. 
\newline
We can adopt a robust approach
considering firstly agent $i$ evolution with no disturbances and adding such effect later. The closed form of the dynamic evolution is given by
\begin{align}
x_i(t)=&A_{ii}^t x_i(0)+\sum_{j=1,j\neq i}^N
\sum_{\tau=0}^{t-1}A_{ii}^{t-\tau-1}A_{ij} x_j(\tau)+ \nonumber\\
& \sum_{j=1}^N
\sum_{\tau=0}^{t-1} A_{ii}^{t-\tau-1}B_{ij} u_j(\tau)+
\sum_{\tau=0}^{t-1} A_{ii}^{t-\tau-1}K_i. \label{eq:linear_evolution}
\end{align} 
\newline
Calling $z_i$ the vector of unknowns such  that $z_i\in \Bbb{R}^{|\mathcal{B}_i^H|}=\mathcal{P}_{\mathcal{B}_i^H}^{\bar{\mathcal{B}}^H}(\mathbf{z})$, with $\mathbf{z}$ provided in Section \ref{subsec:axis_sets_distributed_reachability}, \eqref{eq:linear_evolution} can be rewritten as 
\begin{equation}\label{eq:matrix_F_i_i_system_no_disturbances}
F_i^i z_i=f_i,
\end{equation}
via suitably rearranging the terms in \eqref{eq:linear_evolution} and introducing appropriate null entries for $F_i^i$.
\newline
Now, let us consider matrix
\begin{equation*}
\mathbf{A}=\left[
\begin{array}{ccc}
A_{11} &\dots & A_{1N}\\
\vdots &\vdots & \vdots\\
A_{N1} &\dots & A_{NN}
\end{array}\right].
\end{equation*}
\newline
Since all $\mathcal{L}_{il}(\cdot)\leq 0$ codify linear inequalities and $X_i, U_i$ are polytopes, all the constraints can be rewritten via a set of linear inequality constraints. 
For this reason, the  constraints on the unknown $z_i$ can be written as
\begin{equation}\label{eq:matrix_G_i_i_system_no_disturbances}
G_i^i z_i \leq g_i. 
\end{equation}
\newline
Notice that the solutions $z_i$ satisfying at the same time equality  \eqref{eq:matrix_F_i_i_system_no_disturbances} and inequality \eqref{eq:matrix_G_i_i_system_no_disturbances} are the solutions of the local reachability system \eqref{eq:system_of_constraints_i} for the specific case of affine linear systems on polytopes with no disturbance.
To cope with the disturbances in a robust way, we exploit the superposition property of linear systems.
To do so, let us introduce the block diagonal matrix $\mathbf{E}=\mathrm{diag}\left(E_{1}, \dots, E_{N} \right)$
and matrix 
\begin{equation*}
L_i=[O_{n_1} \dots O_{n_{i-1}} I_{n_i} O_{n_{i+1}} \dots O_{n_N}], 
\end{equation*}
where $O_{n_j}$ and $I_{n_i}$ are the null square matrix and the identity matrix of dimensions $n_j$ and $n_i$, respectively. 
The disturbance action on state $j$ at time $t$, say $\chi_j(t)$ is due to the combined effect of all the disturbances acting on the networked system up to time instant $t-1$ and can be computed, via considering the explicit solution of the overall linear systems, with the expression
\begin{equation}\label{eq:vector_x_j^d}
\chi_j(t)=L_j\sum_{\tau=0}^{t-2}\mathbf{A}^{t-\tau-2}\mathbf{E}d(\tau),
\end{equation}
with $d(t)=[d_1^T(t),d_2^T(t),\dots,d_N^T(t)]^T$ disturbance stack vector.
\newline
Thanks to the above expression, it is possible to consider the vector $\zeta_i$ containing all the $\chi_j(t)$, for $j\in\mathcal{M}_i$ and $t=0,\dots, H$. More precisely, introducing $\mathbf{d}=[d^T(1),d^T(2),\dots,d^T(H)]^T$, vector
$\zeta_i$ assumes expression
\begin{equation*}
\zeta_i=\mathbf{L}_i \mathbf{d},
\end{equation*}
where $\mathbf{L}_i$ is computed\footnote{We do not report the expressions of $F_i^i$ and $\mathbf{L}_i$, since the derivation of the linear affine system case for the theory developed in this paper will be presented in a dedicated contribution.} evaluating the \eqref{eq:vector_x_j^d} according to the specific $j\in\mathcal{M}_i$ and $t$ and placing it in the right block row of $\mathbf{L}_i$ according to the position of $x_j(t)$ in the stack $z_i$. All the rows of $\mathbf{L}_i$ related to the positions of inputs $u_j(t)$ are left null. 
\newline
When including the disturbances over the horizon $t=0, \dots, H$, inequality  \eqref{eq:matrix_G_i_i_system_no_disturbances} is modified as
\begin{equation*}
G_i^i(z_i+\zeta_i)\leq g_i, 
\end{equation*}
where, via simple manipulations we obtain
\begin{equation}\label{eq:matrix_G_i_i_system_disturbances_manipulated}
G_i^iz_i\leq g_i-G_i^i\mathbf{L}_i \mathbf{d}. 
\end{equation}
\newline
Notice that the above expression has to be guaranteed with respect to any admissible disturbance sequence. To do so, we consider the robust solution for system \eqref{eq:matrix_G_i_i_system_disturbances_manipulated}. Calling $n_{g_i}$ the dimension of vector $g_i$, we define 
\begin{equation}\label{eq:delta_i_j}
\delta_i^{(j)}:=\max_{\mathbf{d}\in (D_1\times \dots \times D_N)^H} \nu_j^TG_i^{i}\mathbf{L}_i \mathbf{d},
\end{equation}
and vector $\delta_i:=\left(\delta_i^{(1)}, \delta_i^{(2)}, \dots, \delta_i^{(n_{g_i})}\right)^T$, where $\nu_j\in\Bbb{R}^{n_{g_i}}$ is a vector of a unitary entry at position $j$ and null entries otherwise. Notice that \eqref{eq:delta_i_j} is a LP optimization problem. 
\newline
The robust solution (i.e., the solution feasible with respect to any disturbance sequence) of system \eqref{eq:matrix_G_i_i_system_disturbances_manipulated} is so obtained by solving system
\begin{equation}\label{eq:matrix_G_i_i_system_disturbances_robust}
G_i^iz_i\leq g_i-\delta_i. 
\end{equation}
\newline
It is so easy to notice that, in order to solve the reachability problem in a distributed way, it suffices to iterate equation \eqref{eq:distributed_extrusion_generated_set_no_disturbances} considering $S_{kh,i}(0)$ the solution of the system given by the pair \eqref{eq:matrix_F_i_i_system_no_disturbances},\eqref{eq:matrix_G_i_i_system_disturbances_robust}. 
\newline
Notice also that, when dealing with linear affine systems and linear inequality constraints, the operation involved in equation \eqref{eq:distributed_extrusion_generated_set_no_disturbances} can be easily computed by considering the fact that each neighbour sends to node $i$ its set of constraints. After that, node $i$ will cope with such constraints for a larger vector of unknown variables, say~ $z_{\mathcal{M}_i}\in\Bbb{R}^{|\mathcal{B}_{\mathcal{M}_i}|}$, via suitably including null terms into the linear constraint formulas. Similarly, matrices $G_{i}^i$ and $F_i^i$ and vectors $f_i$, $g_i$ and $\delta_i$ are suitably increased in their dimension via appropriately adding null columns and rows. The whole process  will require iteratively solving only LP problems of size $|\mathcal{B}_{\mathcal{M}_i}^H|$. Also, since the set obtained at each iteration are polytopes (a numerical example on polytopes has been provided in Example \ref{exa:extrusion_generated_polyhedron}), a termination criteria can be considered when the volume of the difference set between two successive iterations is below a given precision threshold. 
\newline
Notice that such results can be then exploited for computing in a distributed way the invariant set of the overall system. Such a problem, extremely relevant in automatic control, has been studied under different more restrictive hypotheses, such as linear autonomous interconnected dynamics in \cite{rake:11}.
Notice also that, at each agent level, the local reachability solution can be directly used as local constraints of a distributed MPC whose aim is the one of picking an (optimal) solution for controlling the overall system.    
Details will be provided in a dedicated paper.

\section{Conclusions}\label{sec:conclusions}
In this work we address the problem of backward distributed reachability for networked nonlinear systems with coupled dynamics and constraints.

We approach the problem from a theoretical perspective. Specifically, we consider a nonlinear constrained large scale system arising from the interconnection of nonlinear agents possibly coupled in their states and inputs though the dynamic function and/or through nonlinear constraints.
No specific hypothesis is required on the dynamics, the constraints or the interconnection pattern.
\newline
For the overall system (here termed as {\em centralized}), we formulate a backward reachability problem and we study the centralized vs distributed relation. Specifically, considering the problem of controlling a large scale system to a target region under suitable control sequences of a given length $H$, we demonstrate through the proposed theory that it is always possible to cast the problem into a distributed equivalent form.
To this aim, suitably mathematical concepts and operators are introduced, with their related formalism. 
\newline
As a result, at each agent level, a reduced size backward reachability problem is derived on a corresponding reduced set of variables depending on the local system. Via cooperatively solving such local problems with information exchange of each agent with its neighbourhood, the centralized problem solution can be obtained with no loss of information. 
\newline
The convergence of the approach is proven, showing that each subsystem is able to compute only the ``portion'' of the overall solution it cares, i.e., the projection of the solution of the overall centralized problem onto a properly defined subspace specific to each agent. 
Furthermore, such portion is the minimal distributed information able to reconstruct the overall centralized solution with no information losses.  
\newline
Under these conditions, indeed, the distributed reachability problem does neither introduce conservativeness nor approximation.
\newline
In addition to this, centralized operations do not need to be performed nor it is needed to reconstruct centrally
the overall solution. The latter is uniquely determined by the computed local distributed projections and satisfies all the dynamic and static constraints. 
In this regard, the proposed approach casts the original centralized reachability problem into a distributed equivalent one. 
\newline
The results presented in this work are purely theoretical and completely agnostic to how the reachability subproblems and set operations can be implemented in practice or which analytical/numerical technique is used for the local problems. 
\newline
A discussion on the applicability of the proposed theory for future research lines is provided in the paper. In particular, some details are given for the linear affine systems case on polytopes. 
It is highlighted how the theory can be conveniently exploited to study backward reachability or computing the control invariant set for any arbitrarily large state space system, with no restrictive hypotheses on the dynamics and the couplings. 
\newline
Casting the theory for some specific class of problems requires, however, dedicated amount of work and its is beyond the scope of the paper.  
Therefore, further details on what reported in the applicability discussion section will be the subject of a dedicated future paper.


%
\begin{ack}                               
The authors wish to warmly thank Angel Molina Acosta for his support in coding the numerical example, as part of his PhD thesis.  
\end{ack}

\bibliographystyle{plain}        
\bibliography{bibliografiagenerale6_21}           



\clearpage
\appendix 
\section*{APPENDIX A: Examples}
\renewcommand{\thefigure}{A.\arabic{figure}}

\begin{figure}
\centering
\includegraphics[width=1\textwidth]{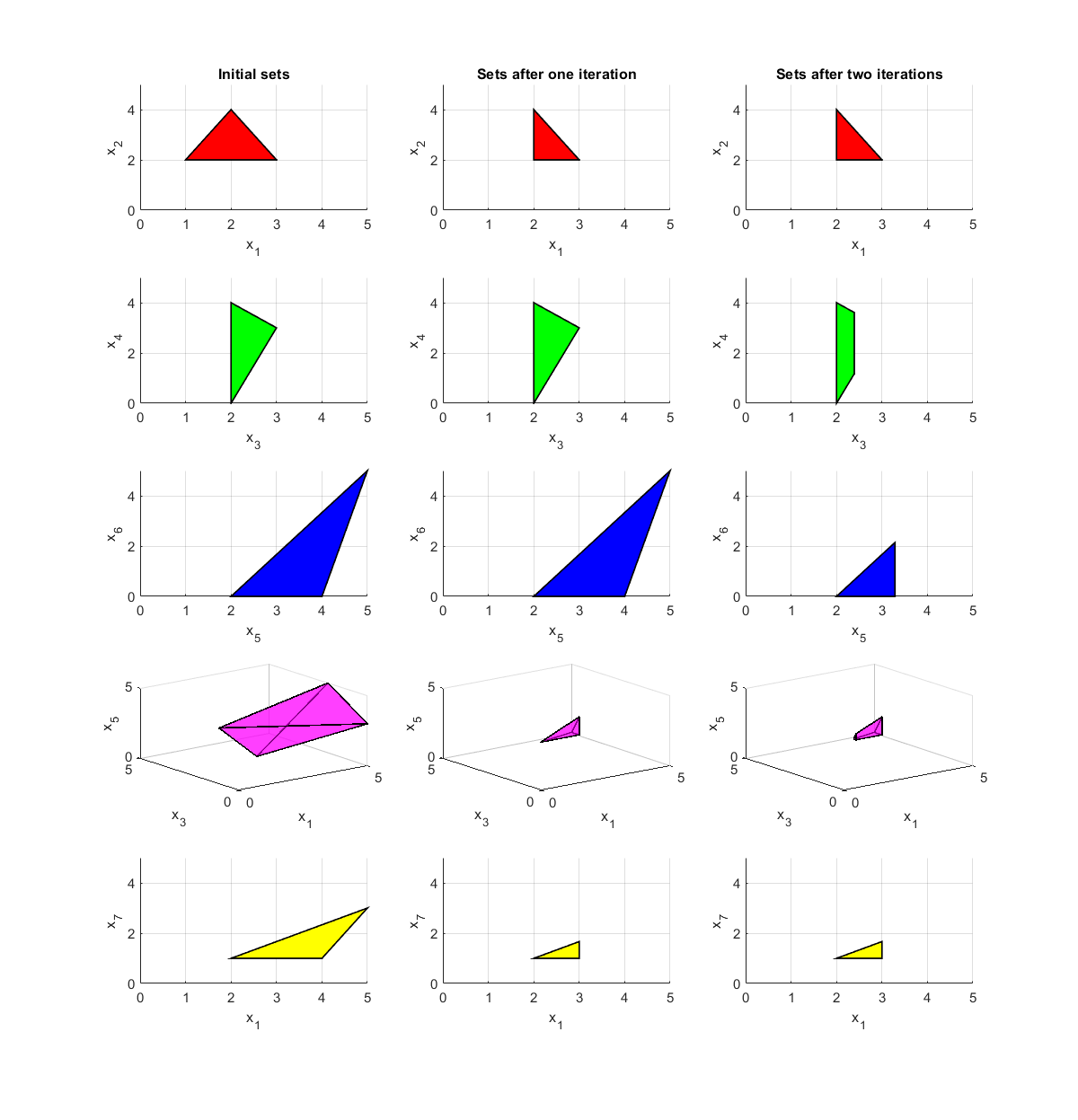}
\caption{Algorithm \ref{alg:distributed_extrusion_generated_set} evolution for Example \ref{exa:extrusion_generated_polyhedron}.}
\label{fig:evolution_polyhedron_example}
\end{figure}

\begin{example}\label{exa:projection_operator}
Consider the axis sets $\mathcal{B}_1=\{3,6,7\}$ and $\mathcal{B}_2=\{1,3,4,6,7\}$. We obviously have $\mathcal{B}_1\subseteq \mathcal{B}_2$, with $|\mathcal{B}_1|=3$ and $|\mathcal{B}_2|=5$. Consider now vector $v=(-4, 6, \pi, 0, 3.2)^T\in\Bbb{R}^5$. We have 
$\mathcal{P}_{\mathcal{B}_1}^{\mathcal{B}_2}(v)=\{(6, 0, 3.2)^T\}$.
\end{example}
\begin{example}\label{exa:extrusion_operator}
Consider the axis sets $\mathcal{B}_1=\{3,5\}$ and $\mathcal{B}_2=\{2,3,4,5,9\}$, with $|\mathcal{B}_1|=2$ and $|\mathcal{B}_2|=5$. Consider now vector $w=(5, -1)^T\in\Bbb{R}^2$. We have 
$\mathcal{E}_{\mathcal{B}_1}^{\mathcal{B}_2}(w)=\{(\alpha_1, 5, \alpha_2,-1, \alpha_3)^T: \, \alpha_1,\alpha_2,\alpha_3 \in \Bbb{R}\}$. 
\end{example}
\begin{example}\label{exa:centralized_computation}
Consider axis sets: $\mathcal{B}_1=\{3,5\}$, $\mathcal{B}_2=\{1,2,3\}$, $\mathcal{B}_3=\{2,5\}$, $\mathcal{B}_4=\{1,4,6\}$, $\mathcal{B}_5=\{4,6\}$.
\newline
Further, consider sets: 
\newline
$S_1=\{(6,7)^T, (-5,-3)^T, (5,-3)^T, (0,0)^T\}$, 
\newline
$S_2=\{(-2,6,5)^T, (1,7,-5)^T, (5,1,0)^T, (5,-1,0)^T\}$,
\newline
$S_3=\{(6,-3)^T, (1,0)^T, (-1,0)^T,(7,6)^T \}$,
\newline
$S_4=\{(-2,0,1)^T, (5,3,-2)^T, (4,3,-2)^T,(1,7,-4)^T \}$,
\newline
$S_5=\{(7,-4)^T, (0,1)^T, (3,-2)^T \}$.
\newline
We have $\bar{\mathcal{B}}=\{1,2,3,4,5,6\}$ and 
\begin{align*}
\bar{S}=&\left\{(-2,6,5,0,-3,1)^T,(5,1,0,3,0,-2)^T, \right.\\
&\left.(5,-1,0,3,0,-2)^T\right\}.
\end{align*}
\end{example}
\begin{example}\label{exa:centralized_computation_from_projection}
Consider the same setting of Example \ref{exa:centralized_computation}. Having computed $\bar{S}$, we have
\newline
$\bar{S}_1=\{(-5,-3)^T, (0,0)^T\}$, 
\newline
$\bar{S}_2=\{(-2,6,5)^T, (5,1,0)^T, (5,-1,0)^T\}$,
\newline
$\bar{S}_3=\{(6,-3)^T, (1,0)^T, (-1,0)^T,\}$,
\newline
$\bar{S}_4=\{(-2,0,1)^T, (5,3,-2)^T \}$,
\newline
$\bar{S}_5=\{(0,1)^T, (3,-2)^T \}$.
\newline
It is possible to see that, computing \eqref{eq:bar_S_centralized_from_projection}, we again obtain $\bar{S}$ as derived in Example \ref{exa:centralized_computation}.
\end{example}
\begin{example}\label{exa:graph}
Consider the axis sets defined in Example~\ref{exa:centralized_computation}. We consider a graph of five nodes such that node $1$ is associated to $\mathcal{B}_1$, node $2$ with $\mathcal{B}_2$ and so on. We therefore have: $\mathcal{M}_1=\{1,2,3\}$, $\mathcal{M}_2=\{1,2,3,4\}$, $\mathcal{M}_3=\{1,2,3\}$, $\mathcal{M}_4=\{2,4,5\}$, $\mathcal{M}_5=\{4,5\}$. 
\end{example}
\begin{example}\label{exa:extrusion_generated_set}
Consider the same setting as in Example \ref{exa:centralized_computation} with the related graph as per Example \ref{exa:graph}. Also, consider as initial sets the same $S_i$ provided in Example \ref{exa:centralized_computation}. Namely
\newline
$S_1(0)=\{(6,7)^T, (-5,-3)^T, (5,-3)^T, (0,0)^T\}$, 
\newline
$S_2(0)=\{(-2,6,5)^T, (1,7,-5)^T, (5,1,0)^T, (5,-1,0)^T\}$,
\newline
$S_3(0)=\{(6,-3)^T, (1,0)^T, (-1,0)^T,(7,6)^T \}$,
\newline
$S_4(0)=\{(-2,0,1)^T, (5,3,-2)^T, (4,3,-2)^T,(1,7,-4)^T \}$,
\newline
$S_5(0)=\{(7,-4)^T, (0,1)^T, (3,-2)^T \}$.
\newline
Iterating Algorithm \ref{alg:distributed_extrusion_generated_set}, we obtain
\newline
$S_1(1)=\{(-5,-3)^T, (0,0)^T\}$, 
\newline
$S_2(1)=\{(-2,6,5)^T, (5,1,0)^T, (5,-1,0)^T\}$,
\newline
$S_3(1)=\{(6,-3)^T, (1,0)^T, (-1,0)^T,\}$,
\newline
$S_4(1)=\{(-2,0,1)^T, (5,3,-2)^T,(1,7,-4)^T \}$,
\newline
$S_5(1)=\{(7,-4)^T,(0,1)^T, (3,-2)^T \}$,
\newline
and
\newline
$S_1(2)=\{(-5,-3)^T, (0,0)^T\}$, 
\newline
$S_2(2)=\{(-2,6,5)^T, (5,1,0)^T, (5,-1,0)^T\}$,
\newline
$S_3(2)=\{(6,-3)^T, (1,0)^T, (-1,0)^T,\}$,
\newline
$S_4(2)=\{(-2,0,1)^T, (5,3,-2)^T \}$,
\newline
$S_5(2)=\{(0,1)^T, (3,-2)^T \}$.
\newline
The latter represent the fixed point of the algorithm, since further iterations provide the same sets as in iteration $\kappa=2$. Notice that, as proved, the above sets are the same as the ones obtained through the centralized operations \eqref{eq:extrusion_generated_set} and \eqref{eq:projection_extrusion_generated_set}. They can be used to reconstruct $\bar{S}$ according to \eqref{eq:bar_S_centralized_from_projection} as showed in Example \ref{exa:centralized_computation_from_projection}. 
\end{example}
\begin{example}\label{exa:extrusion_generated_polyhedron}
Consider a network with $N=5$ nodes with: $\mathcal{B}_1=\{1,2\}$; $\mathcal{B}_2=\{3,4\}$; $\mathcal{B}_3=\{5,6\}$; $\mathcal{B}_4=\{1,3,5\}$; $\mathcal{B}_5=\{2,7\}$.  As per Section \ref{sec:distributed_extrusion_generated_set}, the nodes form the undirected graph $\mathcal{G}_c$ as in Fig. \ref{fig:topology_polyhedron_example}.
Let us consider the simplices sets $S_i$ defined as
{\scriptsize
\begin{align*}
    S_1:=&\left\lbrace\begin{bmatrix}z_1\\z_2\end{bmatrix}\in\mathbb{R}^2:\begin{bmatrix}z_1\\z_2\end{bmatrix}\in\mathbf{Co}\left\lbrace
    \begin{bmatrix}1\\2\end{bmatrix},
    \begin{bmatrix}3\\2\end{bmatrix},
    \begin{bmatrix}2\\4\end{bmatrix}\right\rbrace\right\rbrace,\\
\\
    S_2:=&\left\lbrace\begin{bmatrix}z_3\\z_4\end{bmatrix}\in\mathbb{R}^2:\begin{bmatrix}z_3\\z_4\end{bmatrix}\in\mathbf{Co}\left\lbrace
    \begin{bmatrix}2\\4\end{bmatrix},
    \begin{bmatrix}3\\3\end{bmatrix},
    \begin{bmatrix}2\\0\end{bmatrix}\right\rbrace\right\rbrace,\\
\\
    S_3:=&\left\lbrace\begin{bmatrix}z_5\\z_6\end{bmatrix}\in\mathbb{R}^2:\begin{bmatrix}z_5\\z_6\end{bmatrix}\in\mathbf{Co}\left\lbrace
    \begin{bmatrix}5\\5\end{bmatrix},
    \begin{bmatrix}4\\0\end{bmatrix},
    \begin{bmatrix}2\\0\end{bmatrix}\right\rbrace\right\rbrace,\\
\\
    S_4:=&\left\lbrace\begin{bmatrix}z_1\\z_3\\z_5\end{bmatrix}\in\mathbb{R}^3:\begin{bmatrix}z_1\\z_3\\z_5\end{bmatrix}\in\mathbf{Co}\left\lbrace
    \begin{bmatrix}0\\1\\4\end{bmatrix},
    \begin{bmatrix}3\\3\\0\end{bmatrix},
    \begin{bmatrix}5\\0\\3\end{bmatrix},
    \begin{bmatrix}5\\2\\5\end{bmatrix}\right\rbrace\right\rbrace,\\
\\
    S_5:=&\left\lbrace\begin{bmatrix}z_2\\z_7\end{bmatrix}\in\mathbb{R}^2:\begin{bmatrix}z_2\\z_7\end{bmatrix}\in\mathbf{Co}\left\lbrace
    \begin{bmatrix}2\\1\end{bmatrix},
    \begin{bmatrix}4\\1\end{bmatrix},
    \begin{bmatrix}5\\3\end{bmatrix}\right\rbrace\right\rbrace,\\
\end{align*}
}
where $\mathbf{Co}\{\cdot\}$ denotes the convex hull operation. 
\begin{figure}[h!]
\centering
\includegraphics[width=0.5\textwidth]{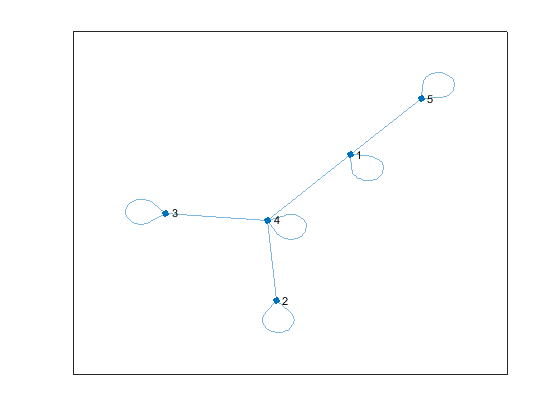}
\caption{Graph topology for Example \ref{exa:extrusion_generated_polyhedron}.}
\label{fig:topology_polyhedron_example}
\end{figure} 
Set $\bar{S}$ computed according to \eqref{eq:extrusion_generated_set} is, therefore, a $7\--$~dimensional polythope. Via computing the centralized construction \eqref{eq:extrusion_generated_set} of such set and the centralized projections $\bar{S}_i$ over the $\mathcal{B}_i$ according to \eqref{eq:projection_extrusion_generated_set} (for the chosen example size it is still possible to perform a centralized calculation, that we use as benchmark), we obtain 
{\scriptsize
\begin{align*}
    \bar{S}_1=&\left\lbrace\begin{bmatrix}z_1\\z_2\end{bmatrix}\in\mathbb{R}^2:\begin{bmatrix}z_1\\z_2\end{bmatrix}\in\mathbf{Co}\left\lbrace
    \begin{bmatrix}2\\2\end{bmatrix},
    \begin{bmatrix}3\\2\end{bmatrix},
    \begin{bmatrix}2\\4\end{bmatrix}\right\rbrace\right\rbrace,\\
\\
    \bar{S}_2:=&\left\lbrace\begin{bmatrix}z_3\\z_4\end{bmatrix}\in\mathbb{R}^2:\begin{bmatrix}z_3\\z_4\end{bmatrix}\in\mathbf{Co}\left\lbrace
    \begin{bmatrix}2\\4\end{bmatrix},
    \begin{bmatrix}2\\0\end{bmatrix},
    \begin{bmatrix}2.39\\1.17\end{bmatrix},
    \begin{bmatrix}2.39\\3.61\end{bmatrix}\right\rbrace\right\rbrace,\\
\\
    \bar{S}_3:=&\left\lbrace\begin{bmatrix}z_5\\z_6\end{bmatrix}\in\mathbb{R}^2:\begin{bmatrix}z_5\\z_6\end{bmatrix}\in\mathbf{Co}\left\lbrace
    \begin{bmatrix}3.29\\0\end{bmatrix},
    \begin{bmatrix}3.29\\2.14\end{bmatrix},
    \begin{bmatrix}2\\0\end{bmatrix}\right\rbrace\right\rbrace,\\
\\
    \bar{S}_4:=&\left\lbrace\begin{bmatrix}z_1\\z_3\\z_5\end{bmatrix}\in\mathbb{R}^3:\begin{bmatrix}z_1\\z_3\\z_5\end{bmatrix}\in\mathbf{Co}\left\lbrace
    \begin{bmatrix}3\\2\\2\end{bmatrix},
    \begin{bmatrix}3\\2.39\\2\end{bmatrix},
    \begin{bmatrix}3\\2\\3.29\end{bmatrix},
    \begin{bmatrix}2\\2\\2.43\end{bmatrix},\right.\right.\\
    &
    \left.\left.\begin{bmatrix}2\\2.13\\2\end{bmatrix},
    \begin{bmatrix}2\\2\\2\end{bmatrix}\right\rbrace\right\rbrace,\\
\\
    \bar{S}_5:=&\left\lbrace\begin{bmatrix}z_2\\z_7\end{bmatrix}\in\mathbb{R}^2:\begin{bmatrix}z_2\\z_7\end{bmatrix}\in\mathbf{Co}\left\lbrace
    \begin{bmatrix}2\\1\end{bmatrix},
    \begin{bmatrix}3\\1\end{bmatrix},
    \begin{bmatrix}3\\1.67\end{bmatrix}\right\rbrace\right\rbrace.\\
\end{align*}
}
When applying the distributed Algorithm \ref{alg:distributed_extrusion_generated_set}, it is possible to observe that it converges in two iterations, which are reported in Fig. \ref{fig:evolution_polyhedron_example}. As it is possible to acknowledge, the projection obtained via the distributed method coincides with the ones computed via centrally reconstructing the set and then projecting it on each of the considered subspaces.
\newline
It is also immediate to notice that the distributed approach can be applied as it is in case the network topology is increased, adding an arbitrary number of further nodes. Conversely, the centralized approach quickly results impractical due to the variables' space growth. 
\end{example}

\onecolumn
\twocolumn

\newpage
\appendix

\section*{APPENDIX B: Theorems proofs}   
\renewcommand{\theequation}{B.\arabic{equation}}
    
\begin{pf}\textbf{(Theorem \ref{thm:bar_S_centralized_from_projection})}
The formula \eqref{eq:bar_S_centralized_from_projection} can be easily proven considering that, by construction, $\bar{S}_i\subseteq S_i$. So, $\bigcap_{i=1}^N \mathcal{E}_{\mathcal{B}_i}^{\bar{\mathcal{B}}}\left(\bar{S}_i \right)\subseteq \bigcap_{i=1}^N \mathcal{E}_{\mathcal{B}_i}^{\bar{\mathcal{B}}}\left(S_i \right)=\bar{S}$. On the other hand, $\bar{S}\subseteq\mathcal{E}_{\mathcal{B}_i}^{\bar{\mathcal{B}}}(\bar{S}_i)$ for all $i=1,\dots,N$ and $\bar{S}\subseteq \bigcap_{i=1}^N \mathcal{E}_{\mathcal{B}_i}^{\bar{\mathcal{B}}}\left(\bar{S}_i \right)$. Hence,
\eqref{eq:bar_S_centralized_from_projection} holds. 

For the second statement, note that, first of all, if $\bar{S}\neq \emptyset$, then $\bar{S}_i\neq \emptyset$ for all $i$. Let us now consider an $\hat{S}_j\subset \bar{S}_j$ and a point $\bar{s}_j\in \bar{S}_j$ such that $\bar{s}_j \notin \hat{S}_j$. Since $\bar{s}_j\in \bar{S}_j$, there exists a point $\bar{s}\in \bar{S}$ such that $\bar{s}_j=\mathcal{P}_{\mathcal{B}_j}^{\bar{\mathcal{B}}}(\bar{s})$. Since the point $\bar{s}\notin \hat{S}$, the strict inclusion $\hat{S}\subset \bar{S}$ holds.
\end{pf}

\begin{pf}\textbf{(Lemma \ref{lem:succession_nested_inclusions})}
To prove {\em i.}, first define 
\newline
$\bar{S}:=\bigcap_{i=1}^N\mathcal{E}_{\mathcal{B}_i}^{\bar{\mathcal{B}}}(S_i(0))$. Using an induction argument, suppose that $\bar{S}=\bigcap_{i=1}^N\mathcal{E}_{\mathcal{B}_i}^{\bar{\mathcal{B}}}(S_i(\kappa))$. This holds at least for $\kappa=0$ by definition. 

It is straightforward to show that $\bar{S}\subseteq \mathcal{E}_{\mathcal{B}_i}^{\bar{\mathcal{B}}}(S_{i}(\kappa+1))$ for all $i=1,\dots, N$. Indeed, for any generic point $\bar{s}\in\bar{S}$, we have that $\bar{s}=\bigcap_{i=1}^N\mathcal{E}_{\mathcal{B}_i}^{\bar{\mathcal{B}}}(\bar{s}_i)$, with $\bar{s}_i=\mathcal{P}_{\mathcal{B}_i}^{\bar{\mathcal{B}}}(\bar{s})$. This implies that, directly from  \eqref{eq:distributed_extrusion_generated_set}, $\bar{s}_i\in S_i(\kappa+1) \,\, \forall i$ and $\bar{S}\subseteq \mathcal{E}_{\mathcal{B}_i}^{\bar{\mathcal{B}}}(S_{i}(\kappa+1))$.  By considering the intersections among all the extrusion sets, we have 
\begin{equation}\label{eq:intersection_left_inclusion_bar_S}
\bar{S}\subseteq \bigcap_{i=1}^N \mathcal{E}_{\mathcal{B}_i}^{\bar{\mathcal{B}}}(S_i(\kappa+1)).
\end{equation}

From \eqref{eq:distributed_extrusion_generated_set}, $S_i(\kappa+1)\subseteq S_i(\kappa) \,\, \forall i$ holds and, by considering the intersections among all the extrusion sets, we can write 

\begin{equation}\label{eq:intersection_right_inclusion_bar_S}
\bigcap_{i=1}^N \mathcal{E}_{\mathcal{B}_i}^{\bar{\mathcal{B}}}(S_i(\kappa+1))\subseteq \bigcap_{i=1}^N \mathcal{E}_{\mathcal{B}_i}^{\bar{\mathcal{B}}}(S_i(\kappa))=\bar{S}. 
\end{equation}

The result in  i. follows from \eqref{eq:intersection_left_inclusion_bar_S} and \eqref{eq:intersection_right_inclusion_bar_S}. 

To prove point ii., we first consider that, from i. and Theorem \ref{thm:bar_S_centralized_from_projection}, $\bar{S}_i\subseteq S_i (\kappa), \,\, \forall i, \, \forall \kappa$. Also, since $S_i(\kappa+1)\subseteq S_i(\kappa), \,\, \forall i, \, \forall \kappa$, then point ii. is proven. 
\end{pf}

\begin{pf}\textbf{(Lemma \ref{lem:fixed_point})}
We will start proving point i. It is straightforward to notice that equation \eqref{eq:distributed_extrusion_generated_set} can be equivalently written as $S_i(\kappa+1)=\mathcal{P}_{\mathcal{B}_i}^{\bar{\mathcal{B}}}\left[
\bigcap_{j\in\mathcal{M}_i}
\mathcal{E}_{\mathcal{B}_j}^{\bar{\mathcal{B}}}\left(S_j(\kappa)\right)
\right]$. Indeed, extending the extrusion operator up to $\bar{\mathcal{B}}$ does not alter the intersection operation since the axis in $\bar{\mathcal{B}}\backslash \mathcal{B}_{\mathcal{M}_i}$ are not shared among the neighbours of node $i$ (entries of points extruded on $\bar{\mathcal{B}}\backslash \mathcal{B}_{\mathcal{M}_i}$ can assume any possible value, thus not affecting the intersection result).
Furthermore, it is also straightforward to recognize that from \eqref{eq:S_i_star}, since $S_i^{*}\subseteq S_i (\kappa)$, we have that $S_i^{*}\cap S_i(\kappa)=S_i^{*}$ for all $\kappa$ and for all $i=1,\dots, N$.

Let us now observe that, for all $\kappa$ and for all $i=1,\dots, N$, it holds that  
\begin{equation}\label{eq:induction_equation_fixed_point_lemma_starting_iteration}
\mathcal{P}_{\mathcal{B}_i}^{\bar{\mathcal{B}}}\left[
\mathcal{E}_{\mathcal{B}_i}^{\bar{\mathcal{B}}}(S_i^{*})\cap
\bigcap_{j\in\mathcal{M}_i-\{i\}}\mathcal{E}_{\mathcal{B}_j}^{\bar{\mathcal{B}}}(S_j(\kappa))
\right]
=S_i^* .
\end{equation}

Indeed, we have

\begin{subequations}
\begin{align*}
&\mathcal{P}_{\mathcal{B}_i}^{\bar{\mathcal{B}}}\left[
\mathcal{E}_{\mathcal{B}_i}^{\bar{\mathcal{B}}}(S_i^{*})
\cap\bigcap_{j\in\mathcal{M}_i-\{i\}}\mathcal{E}_{\mathcal{B}_j}^{\bar{\mathcal{B}}}(S_j(\kappa))
\right]
=\\
&\mathcal{P}_{\mathcal{B}_i}^{\bar{\mathcal{B}}}\left[
\mathcal{E}_{\mathcal{B}_i}^{\bar{\mathcal{B}}}(S_i^{*}\cap S_i(\kappa))
\cap\bigcap_{j\in\mathcal{M}_i-\{i\}}\mathcal{E}_{\mathcal{B}_j}^{\bar{\mathcal{B}}}(S_j(\kappa))
\right]=\\
& \mathcal{P}_{\mathcal{B}_i}^{\bar{\mathcal{B}}}\left[
\bigcap_{j\in\mathcal{M}_i}\mathcal{E}_{\mathcal{B}_j}^{\bar{\mathcal{B}}}(S_j(\kappa))
\right]\cap S_i^*=\\
& S_i(\kappa+1)\cap S_i^*=S_i^*.
\end{align*}
\end{subequations}

Let us now order the elements in $\mathcal{M}_i$ as $\mathcal{M}_i=\{j^1,j^2,\dots,j^{|\mathcal{M}_i|}\}$, with $j^1=i$, and let us consider an index $1\leq h\leq |\mathcal{M}_i|-1$. Using an induction argument, it is possible to prove that if for every $\kappa$ 
\begin{equation}\label{eq:induction_equation_fixed_point_lemma_iteration_l}
\mathcal{P}_{\mathcal{B}_{j^1}}^{\bar{\mathcal{B}}}\left[
\bigcap_{l=1}^h \mathcal{E}_{\mathcal{B}_{j^l}}^{\bar{\mathcal{B}}}(S_{j^l}^*)\cap 
\bigcap_{l=h+1}^{|\mathcal{M}_i|}\mathcal{E}_{\mathcal{B}_{j^l}}^{\bar{\mathcal{B}}}(S_{j^l}(\kappa))
\right]=S_{j^1}^*,
\end{equation}
holds (by equation \eqref{eq:induction_equation_fixed_point_lemma_starting_iteration} it holds at least for $h=1$), then it also holds that
\begin{equation}\label{eq:induction_equation_fixed_point_lemma_iteration_l_plus_one}
\mathcal{P}_{\mathcal{B}_{j^1}}^{\bar{\mathcal{B}}}\left[
\bigcap_{l=1}^{h+1} \mathcal{E}_{\mathcal{B}_{j^l}}^{\bar{\mathcal{B}}}(S_{j^l}^*)\cap 
\bigcap_{l=h+2}^{|\mathcal{M}_i|}\mathcal{E}_{\mathcal{B}_{j^l}}^{\bar{\mathcal{B}}}(S_{j^l}(\kappa))
\right]=S_{j^1}^* . 
\end{equation} 

To prove such implication, we will need two intermediate relations. In view of this, let us first consider a $\kappa^\prime> \kappa$. Equation \eqref{eq:induction_equation_fixed_point_lemma_iteration_l} can be written as  
\begin{align}
&\mathcal{P}_{\mathcal{B}_{j^1}}^{\bar{\mathcal{B}}}\left[
\bigcap_{l=1}^h \mathcal{E}_{\mathcal{B}_{j^l}}^{\bar{\mathcal{B}}}(S_{j^l}^*)
\cap 
\mathcal{E}_{\mathcal{B}_{j^{h+1}}}^{\bar{\mathcal{B}}}(S_{j^{h+1}}(\kappa^\prime))
\cap \right. \nonumber \\
& \left. \bigcap_{l=h+2}^{|\mathcal{M}_i|}\mathcal{E}_{\mathcal{B}_{j^l}}^{\bar{\mathcal{B}}}(S_{j^l}(\kappa))
\right]=S_{j^1}^* . \label{eq:induction_equation_fixed_point_lemma_iteration_l_modified}
\end{align}
This fact follows from
\begin{align*}
S_{j^1}^* =& \mathcal{P}_{\mathcal{B}_{j^1}}^{\bar{\mathcal{B}}}\left[
\bigcap_{l=1}^h \mathcal{E}_{\mathcal{B}_{j^l}}^{\bar{\mathcal{B}}}(S_{j^l}^*)\cap 
\bigcap_{l=h+1}^{|\mathcal{M}_i|}\mathcal{E}_{\mathcal{B}_{j^l}}^{\bar{\mathcal{B}}}(S_{j^l}(\kappa^\prime))
\right]\subseteq \\
 &\mathcal{P}_{\mathcal{B}_{j^1}}^{\bar{\mathcal{B}}}\left[
\bigcap_{l=1}^h \mathcal{E}_{\mathcal{B}_{j^l}}^{\bar{\mathcal{B}}}(S_{j^l}^*)
\cap
\mathcal{E}_{\mathcal{B}_{j^{h+1}}}^{\bar{\mathcal{B}}}(S_{j^{h+1}}(\kappa^\prime)) \cap \right. \\
&\left. 
\bigcap_{l=h+2}^{|\mathcal{M}_i|}\mathcal{E}_{\mathcal{B}_{j^l}}^{\bar{\mathcal{B}}}(S_{j^l}(\kappa))
\right] \subseteq \\
 &\mathcal{P}_{\mathcal{B}_{j^1}}^{\bar{\mathcal{B}}}\left[
\bigcap_{l=1}^h \mathcal{E}_{\mathcal{B}_{j^l}}^{\bar{\mathcal{B}}}(S_{j^l}^*)\cap 
\bigcap_{l=h+1}^{|\mathcal{M}_i|}\mathcal{E}_{\mathcal{B}_{j^l}}^{\bar{\mathcal{B}}}(S_{j^l}(\kappa))
\right]= S_{j^1}^*, 
\end{align*}
where in the first and last member we exploited the equation \eqref{eq:induction_equation_fixed_point_lemma_iteration_l}. 

Consider any $\kappa^\prime>\kappa$ and any $\tilde{S}_{j^{h+1}}\subseteq S_{j^{h+1}}(\kappa)$, the following equation
\begin{align} \label{eq:fixed_point_lemma_time_splitting_equation}
&\mathcal{P}_{\mathcal{B}_{j^1}}^{\bar{\mathcal{B}}} \left\{
\bigcap_{l=1}^h\mathcal{E}_{\mathcal{B}_{j^l}}^{\bar{\mathcal{B}}}(S_{j^l}^*)\cap 
\mathcal{E}_{\mathcal{B}_{j^{h+1}}}^{\bar{\mathcal{B}}}\left[S_{j^{h+1}}(\kappa^\prime)\cap \tilde{S}_{j^{h+1}} \right]\cap \right. \nonumber \\
&\left. \bigcap_{l=h+2}^{|\mathcal{M}_i|}\mathcal{E}_{\mathcal{B}_{j^l}^{\bar{\mathcal{B}}}}(S_{j^l}(\kappa))
\right\}=  \nonumber \\
&\mathcal{P}_{\mathcal{B}_{j^1}}^{\bar{\mathcal{B}}}\left[
\bigcap_{l=1}^h\mathcal{E}_{\mathcal{B}_{j^l}}^{\bar{\mathcal{B}}}(S_{j^l}^*)\cap 
\mathcal{E}_{\mathcal{B}_{j^{h+1}}}^{\bar{\mathcal{B}}}(S_{j^{h+1}}(\kappa^\prime))\cap \right. \nonumber\\
&\left. \bigcap_{l=h+2}^{|\mathcal{M}_i|}\mathcal{E}_{\mathcal{B}_{j^l}^{\bar{\mathcal{B}}}}(S_{j^l}(\kappa))
\right]\cap  \mathcal{P}_{\mathcal{B}_{j^1}}^{\bar{\mathcal{B}}}\left[
\bigcap_{l=1}^h\mathcal{E}_{\mathcal{B}_{j^l}}^{\bar{\mathcal{B}}}(S_{j^l}^*)\cap \right. \nonumber\\
&\left. 
\mathcal{E}_{\mathcal{B}_{j^{h+1}}}^{\bar{\mathcal{B}}}(\tilde{S}_{j^{h+1}})\cap
\bigcap_{l=h+2}^{|\mathcal{M}_i|}\mathcal{E}_{\mathcal{B}_{j^l}^{\bar{\mathcal{B}}}}(S_{j^l}(\kappa))
\right].  \nonumber \\
\end{align} 
can be obtained by observing that 
\begin{align*} 
&\mathcal{P}_{\mathcal{B}_{j^1}}^{\bar{\mathcal{B}}}\left[
\bigcap_{l=1}^h\mathcal{E}_{\mathcal{B}_{j^l}}^{\bar{\mathcal{B}}}(S_{j^l}^*)\cap 
\mathcal{E}_{\mathcal{B}_{j^{h+1}}}^{\bar{\mathcal{B}}}(\tilde{S}_{j^{h+1}})\cap\right. \nonumber \\
&\left.
\bigcap_{l=h+2}^{|\mathcal{M}_i|}\mathcal{E}_{\mathcal{B}_{j^l}^{\bar{\mathcal{B}}}}(S_{j^l}(\kappa))
\right] \subseteq \\
&\mathcal{P}_{\mathcal{B}_{j^1}}^{\bar{\mathcal{B}}}\left[
\bigcap_{l=1}^h\mathcal{E}_{\mathcal{B}_{j^l}}^{\bar{\mathcal{B}}}(S_{j^l}^*)\cap 
\mathcal{E}_{\mathcal{B}_{j^{h+1}}}^{\bar{\mathcal{B}}}(S_{j^{h+1}}(\kappa^{\prime}))\cap
\bigcap_{l=h+2}^{|\mathcal{M}_i|}\mathcal{E}_{\mathcal{B}_{j^l}^{\bar{\mathcal{B}}}}(S_{j^l}(\kappa))
\right],
\end{align*} 
which leads to
\begin{align*} 
&\mathcal{P}_{\mathcal{B}_{j^1}}^{\bar{\mathcal{B}}} \left\{
\bigcap_{l=1}^h\mathcal{E}_{\mathcal{B}_{j^l}}^{\bar{\mathcal{B}}}(S_{j^l}^*)\cap 
\mathcal{E}_{\mathcal{B}_{j^{h+1}}}^{\bar{\mathcal{B}}}\left[S_{j^{h+1}}(\kappa^\prime)\cap \tilde{S}_{j^{h+1}} \right]\cap 
\right. \nonumber \\
&\left. 
\bigcap_{l=h+2}^{|\mathcal{M}_i|}\mathcal{E}_{\mathcal{B}_{j^l}^{\bar{\mathcal{B}}}}(S_{j^l}(\kappa))
\right\}=  \\
&\mathcal{P}_{\mathcal{B}_{j^1}}^{\bar{\mathcal{B}}}\left[
\bigcap_{l=1}^h\mathcal{E}_{\mathcal{B}_{j^l}}^{\bar{\mathcal{B}}}(S_{j^l}^*)\cap 
\mathcal{E}_{\mathcal{B}_{j^{h+1}}}^{\bar{\mathcal{B}}}(\tilde{S}_{j^{h+1}})\cap \right. \\
&\left.
\bigcap_{l=h+2}^{|\mathcal{M}_i|}\mathcal{E}_{\mathcal{B}_{j^l}^{\bar{\mathcal{B}}}}(S_{j^l}(\kappa))
\right]=\\
&\mathcal{P}_{\mathcal{B}_{j^1}}^{\bar{\mathcal{B}}}\left[
\bigcap_{l=1}^h\mathcal{E}_{\mathcal{B}_{j^l}}^{\bar{\mathcal{B}}}(S_{j^l}^*)\cap 
\mathcal{E}_{\mathcal{B}_{j^{h+1}}}^{\bar{\mathcal{B}}}(\tilde{S}_{j^{h+1}})\cap\right. \\
&\left. 
\bigcap_{l=h+2}^{|\mathcal{M}_i|}\mathcal{E}_{\mathcal{B}_{j^l}^{\bar{\mathcal{B}}}}(S_{j^l}(\kappa))
\right]\cap 
\mathcal{P}_{\mathcal{B}_{j^1}}^{\bar{\mathcal{B}}}\left[
\bigcap_{l=1}^h\mathcal{E}_{\mathcal{B}_{j^l}}^{\bar{\mathcal{B}}}(S_{j^l}^*)\cap \right. \\
&\left. 
\mathcal{E}_{\mathcal{B}_{j^{h+1}}}^{\bar{\mathcal{B}}}(S_{j^{h+1}}(\kappa^{\prime}))\cap
\bigcap_{l=h+2}^{|\mathcal{M}_i|}\mathcal{E}_{\mathcal{B}_{j^l}^{\bar{\mathcal{B}}}}(S_{j^l}(\kappa))
\right]. 
\end{align*} 
thus proving the \eqref{eq:fixed_point_lemma_time_splitting_equation}.

Let us now consider the expression
\begin{equation}\label{eq:induction_equation_fixed_point_lemma_iteration_l_plus_one_modified}
\mathcal{P}_{\mathcal{B}_{j^1}}^{\bar{\mathcal{B}}} \left\{
\bigcap_{l=1}^h\mathcal{E}_{\mathcal{B}_{j^l}}^{\bar{\mathcal{B}}}(S_{j^l}^*)\cap 
\mathcal{E}_{\mathcal{B}_{j^{h+1}}}^{\bar{\mathcal{B}}} (S_{j^{h+1}}^*)
\cap 
\bigcap_{l=h+2}^{|\mathcal{M}_i|}\mathcal{E}_{\mathcal{B}_{j^l}^{\bar{\mathcal{B}}}}(S_{j^l}(\kappa))
\right\}.
\end{equation}

By explicitly writing $S_{j^{h+1}}^*$ according to \eqref{eq:S_i_star}, we have
\begin{align*}
&\mathcal{P}_{\mathcal{B}_{j^1}}^{\bar{\mathcal{B}}} \left\{
\bigcap_{l=1}^h\mathcal{E}_{\mathcal{B}_{j^l}}^{\bar{\mathcal{B}}}(S_{j^l}^*)\cap 
\mathcal{E}_{\mathcal{B}_{j^{h+1}}}^{\bar{\mathcal{B}}} (S_{j^{h+1}}^*)
\cap 
\bigcap_{l=h+2}^{|\mathcal{M}_i|}\mathcal{E}_{\mathcal{B}_{j^l}^{\bar{\mathcal{B}}}}(S_{j^l}(\kappa))
\right\}=\\
&\mathcal{P}_{\mathcal{B}_{j^1}}^{\bar{\mathcal{B}}} \left\{
\bigcap_{l=1}^h\mathcal{E}_{\mathcal{B}_{j^l}}^{\bar{\mathcal{B}}}(S_{j^l}^*)\cap 
\mathcal{E}_{\mathcal{B}_{j^{h+1}}}^{\bar{\mathcal{B}}} \left[
\bigcap_{\rho=0}^{+\infty} S_{j^{h+1}}(\rho)
\right]
\cap \right. \\
&\left. 
\bigcap_{l=h+2}^{|\mathcal{M}_i|}\mathcal{E}_{\mathcal{B}_{j^l}^{\bar{\mathcal{B}}}}(S_{j^l}(\kappa))
\right\}=\\
&\mathcal{P}_{\mathcal{B}_{j^1}}^{\bar{\mathcal{B}}} \left\{
\bigcap_{l=1}^h\mathcal{E}_{\mathcal{B}_{j^l}}^{\bar{\mathcal{B}}}(S_{j^l}^*)\cap 
\mathcal{E}_{\mathcal{B}_{j^{h+1}}}^{\bar{\mathcal{B}}} \left[
\bigcap_{\rho=\kappa}^{+\infty} S_{j^{h+1}}(\rho)
\right]
\cap \right. \\
&\left. 
\bigcap_{l=h+2}^{|\mathcal{M}_i|}\mathcal{E}_{\mathcal{B}_{j^l}^{\bar{\mathcal{B}}}}(S_{j^l}(\kappa))
\right\},
\end{align*}
where we considered that, due to Lemma \ref{lem:succession_nested_inclusions},
for any $\kappa\geq 0$, $\bigcap_{\rho=0}^{+\infty} S_{j^{h+1}}(\rho)=\bigcap_{\rho=\kappa}^{+\infty} S_{j^{h+1}}(\rho)$. 

From the above formula and by applying the \eqref{eq:fixed_point_lemma_time_splitting_equation} at any step $\kappa$ we have
\begin{align*}
&\mathcal{P}_{\mathcal{B}_{j^1}}^{\bar{\mathcal{B}}} \left\{
\bigcap_{l=1}^h\mathcal{E}_{\mathcal{B}_{j^l}}^{\bar{\mathcal{B}}}(S_{j^l}^*)\cap 
\mathcal{E}_{\mathcal{B}_{j^{h+1}}}^{\bar{\mathcal{B}}} \left[
\bigcap_{\rho=\kappa}^{+\infty} S_{j^{h+1}}(\rho)
\right]
\cap \right. \\
&\left. 
\bigcap_{l=h+2}^{|\mathcal{M}_i|}\mathcal{E}_{\mathcal{B}_{j^l}^{\bar{\mathcal{B}}}}(S_{j^l}(\kappa))
\right\}=\\
&\bigcap_{\rho=\kappa}^{+\infty}\mathcal{P}_{\mathcal{B}_{j^1}}^{\bar{\mathcal{B}}} \left\{
\bigcap_{l=1}^h\mathcal{E}_{\mathcal{B}_{j^l}}^{\bar{\mathcal{B}}}(S_{j^l}^*)\cap 
\mathcal{E}_{\mathcal{B}_{j^{h+1}}}^{\bar{\mathcal{B}}} (
 S_{j^{h+1}}(\rho)
)
\cap \right. \\
&\left. 
\bigcap_{l=h+2}^{|\mathcal{M}_i|}\mathcal{E}_{\mathcal{B}_{j^l}^{\bar{\mathcal{B}}}}(S_{j^l}(\kappa))
\right\}
\end{align*}
Finally, considering \eqref{eq:induction_equation_fixed_point_lemma_iteration_l_modified}, the above formula can be rewritten as
\begin{align*}
&\mathcal{P}_{\mathcal{B}_{j^1}}^{\bar{\mathcal{B}}} \left\{
\bigcap_{l=1}^h\mathcal{E}_{\mathcal{B}_{j^l}}^{\bar{\mathcal{B}}}(S_{j^l}^*)\cap 
\mathcal{E}_{\mathcal{B}_{j^{h+1}}}^{\bar{\mathcal{B}}} \left[
\bigcap_{\rho=\kappa}^{+\infty} S_{j^{h+1}}(\rho)
\right]
\cap \right. \\
&\left. 
\bigcap_{l=h+2}^{|\mathcal{M}_i|}\mathcal{E}_{\mathcal{B}_{j^l}^{\bar{\mathcal{B}}}}(S_{j^l}(\kappa))
\right\}=
\bigcap_{\rho=\kappa}^{+\infty}S_{j^1}^*=S_{j^1}^*,
\end{align*}
thus allowing to state the equality \eqref{eq:induction_equation_fixed_point_lemma_iteration_l_plus_one_modified}$\,=S_{j^1}^*$. Since \eqref{eq:induction_equation_fixed_point_lemma_iteration_l_plus_one_modified} is equal to \eqref{eq:induction_equation_fixed_point_lemma_iteration_l_plus_one}, then the implication \eqref{eq:induction_equation_fixed_point_lemma_iteration_l}$ \Rightarrow$\eqref{eq:induction_equation_fixed_point_lemma_iteration_l_plus_one} is proven. 

The proof of point i. is then obtained by iterating the equation \eqref{eq:induction_equation_fixed_point_lemma_iteration_l}$\--$\eqref{eq:induction_equation_fixed_point_lemma_iteration_l_plus_one} up to $h=|\mathcal{M}_i|-1$.

To prove point ii. let us first consider that, from \eqref{eq:S_i_star} and Lemma \ref{lem:succession_nested_inclusions}, we have $\bar{S}_i\subseteq S_i^*$. Also, since $S_i^*\subseteq S_i(\kappa)$ for any $\kappa$, we can write $\bar{S}_i\subseteq S_i^*\subseteq S_i(\kappa)$. By applying the extrusion operator and then intersecting for all $i$, we obtain
\begin{equation*}
\bigcap_{i=1}^N\mathcal{E}_{\mathcal{B}_i}^{\bar{\mathcal{B}}}(\bar{S}_i)\subseteq
\bigcap_{i=1}^N\mathcal{E}_{\mathcal{B}_i}^{\bar{\mathcal{B}}}(S_i^*)\subseteq
\bigcap_{i=1}^N\mathcal{E}_{\mathcal{B}_i}^{\bar{\mathcal{B}}}(S(\kappa)).
\end{equation*} 

From the above relation and from Lemma \ref{lem:succession_nested_inclusions}, we obtain $\bar{S}\subseteq S^* \subseteq \bar{S}$ and so $S^*=\bar{S}$. This concludes the proof.
\end{pf}

\setcounter{lemma}{0}
\renewcommand{\thelemma}{B.\arabic{lemma}}
\begin{lemma}\label{lem:extracted_set}
Let us consider any subset $\tilde{S}_i^*\subseteq S_i^*$, with $S_i^*$ fixed point of Algorithm \ref{alg:distributed_extrusion_generated_set} according to  Lemma \ref{lem:fixed_point}. If $\tilde{S}_i^*\neq \emptyset$, then $\forall j\in\mathcal{M}_i$ there exists a  subset $\tilde{S}_j^*\subseteq S_j^*$ such that $\tilde{S}_j^*\neq \emptyset$ and $\mathcal{P}_{\mathcal{B}_i\cap\mathcal{B}_j}^{\bar{\mathcal{B}}}\left[\mathcal{E}_{\mathcal{B}_j}^{\bar{\mathcal{B}}}(\tilde{S}_j^*)\right]\subseteq \mathcal{P}_{\mathcal{B}_i\cap\mathcal{B}_j}^{\bar{\mathcal{B}}}\left[\mathcal{E}_{\mathcal{B}_i}^{\bar{\mathcal{B}}}(\tilde{S}_i^*)\right]$.

We will say that the nonempty set $\tilde{S}_i^*$  {\em induces} the subset $\tilde{S}_j^*$ on $S_j^*$.
\end{lemma}

\begin{pf}\textbf{(Lemma \ref{lem:extracted_set})}
From Lemma \ref{lem:fixed_point} point i. (taking into account the equivalence in considering $\bar{\mathcal{B}}$ instead of $\mathcal{B}_{\mathcal{M}_i}$) we have 
\begin{equation}
S_i^{*}= \mathcal{P}_{\mathcal{B}_i}^{\bar{\mathcal{B}}}\left[\bigcap_{j\in \mathcal{M}_i}\mathcal{E}_{\mathcal{B}_j}^{\bar{\mathcal{B}}} (S_j^{*})\right]. 
\end{equation}

Then, by considering $\tilde{S}_i^*$, from the above expression and after simple manipulations we obtain
\begin{equation}\label{eq:intermediate_eq_extracted_lemma}
\tilde{S}_i^*=\mathcal{P}_{\mathcal{B}_i}^{\bar{\mathcal{B}}}\left[
\mathcal{E}_{\mathcal{B}_i}^{\bar{\mathcal{B}}}(\tilde{S}_i^*)\cap
\bigcap_{j\in \mathcal{M}_i-\{i\}}\mathcal{E}_{\mathcal{B}_j}^{\bar{\mathcal{B}}} (S_j^{*})\right]. 
\end{equation}

Indeed, the formula \eqref{eq:intermediate_eq_extracted_lemma} has been obtained by considering the following steps
\begin{align*}
& \mathcal{P}_{\mathcal{B}_i}^{\bar{\mathcal{B}}}\left[
\mathcal{E}_{\mathcal{B}_i}^{\bar{\mathcal{B}}}(\tilde{S}_i^*)
\cap
\bigcap_{j\in \mathcal{M}_i-\{i\}}\mathcal{E}_{\mathcal{B}_j}^{\bar{\mathcal{B}}} (S_j^{*})\right]=\\
& \mathcal{P}_{\mathcal{B}_i}^{\bar{\mathcal{B}}}\left[
\mathcal{E}_{\mathcal{B}_i}^{\bar{\mathcal{B}}}(\tilde{S}_i^*\cap S_i^*)
\cap
\bigcap_{j\in \mathcal{M}_i-\{i\}}\mathcal{E}_{\mathcal{B}_j}^{\bar{\mathcal{B}}} (S_j^{*})\right]=\\
& \mathcal{P}_{\mathcal{B}_i}^{\bar{\mathcal{B}}}\left[
\mathcal{E}_{\mathcal{B}_i}^{\bar{\mathcal{B}}}(S_i^*)
\cap
\bigcap_{j\in \mathcal{M}_i-\{i\}}\mathcal{E}_{\mathcal{B}_j}^{\bar{\mathcal{B}}} (S_j^{*})\right]\cap \tilde{S}_i^*=\\
& S_i^*\cap \tilde{S}_i^*=\tilde{S}_i^*.
\end{align*}

From \eqref{eq:intermediate_eq_extracted_lemma} it is clear that, if $\tilde{S}_i^*\neq \emptyset$, then $\mathcal{E}_{\mathcal{B}_i}^{\bar{\mathcal{B}}}(\tilde{S}_i^*)\cap\mathcal{E}_{\mathcal{B}_j}^{\bar{\mathcal{B}}}(S_j^*)\neq \emptyset,\,\,\forall j\in\mathcal{M}_i$.
This implies that points $s\in\mathcal{E}_{\mathcal{B}_j}^{\bar{\mathcal{B}}}(S_j^*)$ exist such that  $s\in\mathcal{E}_{\mathcal{B}_i}^{\bar{\mathcal{B}}}(\tilde{S}_i^*)$. Hence, 
\begin{equation*}
\mathcal{P}_{\mathcal{B}_i\cap\mathcal{B}_j}^{\bar{\mathcal{B}}}(s)
\in
\mathcal{P}_{\mathcal{B}_i\cap\mathcal{B}_j}^{\bar{\mathcal{B}}}\left[\mathcal{E}_{\mathcal{B}_i}^{\bar{\mathcal{B}}}(\tilde{S}_i^*)\right]
\end{equation*}
and 
\begin{equation*}
\mathcal{P}_{\mathcal{B}_i\cap\mathcal{B}_j}^{\bar{\mathcal{B}}}(s)
\in
\mathcal{P}_{\mathcal{B}_i\cap\mathcal{B}_j}^{\bar{\mathcal{B}}}\left[\mathcal{E}_{\mathcal{B}_j}^{\bar{\mathcal{B}}}(S_j^*)\right].
\end{equation*}

By denoting with $\tilde{S}_j^*\subseteq S_j^*$ the projection of all the points $s$ towards $\mathcal{B}_j$ the lemma is proven.
\end{pf}

\setcounter{theorem}{0}
\renewcommand{\thetheorem}{B.\arabic{theorem}}
\begin{theorem}\label{thm:iterative_extracted_sets}
Let us consider any nonempty subset $\tilde{S}_i^*\subseteq S_i^*$, with $i=1,\dots,N$ and with $S_i^*$ fixed point of Algorithm \ref{alg:distributed_extrusion_generated_set} according to  Lemma \ref{lem:fixed_point}. 
Then: 
\begin{itemize}
\item[1)] $\tilde{S}_i^*$ induces nonempty subsets $\tilde{S}_j^*\subseteq S_j^*$ for all $j=1, \dots, N$ such that $\bigcap_{j=1}^N \mathcal{E}_{\mathcal{B}_j}^{\bar{\mathcal{B}}}(\tilde{S}_j^*)=\tilde{S}^*\neq\emptyset$ and $\mathcal{P}_{\mathcal{B}_i}^{\bar{\mathcal{B}}}(\tilde{S}^*)=\tilde{S}_i^*$.

\item[2)] Furthermore, let us have $M\leq N$ nodes ordered according to any arbitrary sequence $j^1,j^2,\dots,j^M \in R$ and axis sets $\tilde{\mathcal{B}}_{j^1}\subseteq \mathcal{B}_{j^1},\dots, \tilde{\mathcal{B}}_{j^M}\subseteq \mathcal{B}_{j^M}$. Given the points $\tilde{z}_{j^1}\in \Bbb{R}^{|\tilde{\mathcal{B}}_{j^1}|}, \dots, \tilde{z}_{j^M}\in \Bbb{R}^{|\tilde{\mathcal{B}}_{j^M}|}$, let us 
define the set of points 
\begin{align*}
Z:=& \left\{z\in\bigcap_{j\in\{j^1,\dots, j^M\}}\mathcal{E}_{\mathcal{B}_j}^{\bar{\mathcal{B}}}(\tilde{S}_j^*): \mathcal{P}_{\tilde{\mathcal{B}}_j}^{\bar{\mathcal{B}}}(z) = \tilde{z}_{j}, \right. \\
& \left. \quad\forall j=j^1,\dots, j^M \right\}.
\end{align*} 
 
If $Z\neq \emptyset$, then the two following sentences are equivalent:
\begin{itemize}
\item[(i)] there exists a point $w_1\in\tilde{S}_h^*$, with $h=1,\dots, N$, such that $\mathcal{P}_{\tilde{\mathcal{B}}_j\cap\mathcal{B}_h}^{\mathcal{B}_h}(w_1)=\mathcal{P}_{\tilde{\mathcal{B}}_j\cap\mathcal{B}_h}^{\tilde{\mathcal{B}}_j}(\tilde{z}_j)$, for all $j=j^1,\dots,j^M$. \\
\item[(ii)] there exists a point $w_2\in\tilde{S}_h^*$, with $h=1,\dots, N$, such that $Z\cap \mathcal{E}_{\mathcal{B}_h}^{\bar{\mathcal{B}}}(w_2)\neq \emptyset$. 
\end{itemize}
\end{itemize}
\end{theorem}

\begin{pf}\textbf{(Theorem \ref{thm:iterative_extracted_sets})}
Let us reorder the nodes of the graph as $j^1,j^2, \dots, j^N$. We show that a generic subset $\emptyset\neq\tilde{S}_{j^1}^*\subseteq S_{j^1}^*$ induces nonempty subsets $\tilde{S}_j^*$ on the neighbours of nodes $j^1$ and $j^2$, i.e, $\mathcal{M}_{j^1}$ and $\mathcal{M}_{j^2}$, such as $\bigcap_{j\in(\mathcal{M}_{j^1}\cup \mathcal{M}_{j^2})}\mathcal{E}_{\mathcal{B}_j}^{\bar{\mathcal{B}}}(\tilde{S}_j^*)\neq \emptyset$.
To do so, let us consider the iterative constructive process detailed next. 

As first step of such process, we define
\begin{equation}\label{eq:mathsf_E_j2_prime}
\mathsf{E}_{j^2}^\prime :=\emptyset,  
\end{equation}
and 

\begin{subequations}
\begin{align}
&\eta_{j^1} :=\mathcal{M}_{j^1}, \label{eq:eta_j1}\\
&\mathsf{B}_{j^2} :=\mathcal{M}_{j^2}\cap\eta_{j^1}, \label{eq:mathsf_B_j2}\\
&\mathsf{E}_{j^2} :=\left\{j\in (\mathcal{M}_{j^2}\backslash\mathsf{B}_{j^2}): \mathcal{M}_j\cap(\eta_{j^1}\backslash\mathsf{B}_{j^2})\neq \emptyset \right\},  \label{eq:mathsf_E_j2}\\
&\mathsf{C}_{j^2} :=(\mathcal{M}_{j^2}\backslash(\mathsf{B}_{j^2}\cup \mathsf{E}_{j^2}))\cup\mathsf{E}_{j^2}^\prime, \label{eq:mathsf_C_j2}\\
&\mathsf{D}_{j^1}  :=\left\{j\in (\eta_{j^1}\backslash\mathsf{B}_{j^2}): \mathcal{M}_j\cap \mathsf{E}_{j^2}\neq \emptyset \right\},\label{eq:mathsf_D_j1}\\
&\mathsf{A}_{j^1}  := \eta_{j^1}\backslash(\mathsf{B}_{j^2}\cup \mathsf{D}_{j^1}). \label{eq:mathsf_A_j1}
\end{align}
\end{subequations}

Note that some of the above sets may be empty. For example, in case $j^1$ and $j^2$ neither are not neighbour nor share a common third node in their neighbourhood, $\mathsf{B}_{j^2}$ results empty. 

It is convenient to explicitly enumerate the elements of $\mathsf{E}_{j^2}$ in \eqref{eq:mathsf_E_j2}, i.e., 
\begin{equation}\label{eq:mathsf_E_j2_enumerated_elements}
\mathsf{E}_{j^2}=\left\{e_{j^2}^1,e_{j^2}^2, \dots, e_{j^2}^{|\mathsf{E}_{j^2}|}\right\}.
\end{equation}

We now define $\mathcal{B}_{\eta_{j^1}}=\bigcup_{j\in\eta_{j^1}}\mathcal{B}_j$ and, similarly, $\mathcal{B}_{\mathsf{A}_{j^1}}$, $\mathcal{B}_{\mathsf{D}_{j^1}}$, $\mathcal{B}_{\mathsf{B}_{j^2}}$, $\mathcal{B}_{\mathsf{C}_{j^2}}$. Also, we define $\hat{\mathcal{B}}=\mathcal{B}_{\eta_{j^1}}\cup \mathcal{B}_{\mathsf{B}_{j^2}}\cup\mathcal{B}_{\mathsf{C}_{j^2}}\cup\mathcal{B}_{e_{j^2}^1}$. 

By relying on the above definitions, we consider the following three-equations system
\begin{subequations}
\begin{align}
&S_{\eta_{j^1}}(\kappa+1)=\mathcal{P}_{\mathcal{B}_{\eta_{j^1}}}^{\hat{\mathcal{B}}}\left[
\mathcal{E}_{\mathcal{B}_{\eta_{j^1}}}^{\hat{\mathcal{B}}}(S_{\eta_{j^1}}(\kappa))\cap \right. \nonumber\\
&\left. 
\mathcal{E}_{{\mathcal{B}_{\mathsf{B}_{j^2}}\cup \mathcal{B}_{\mathsf{C}_{j^2}}}}^{\hat{\mathcal{B}}}(S_{
\mathsf{B}_{j^2}\mathsf{C}_{j^2}}(\kappa))\cap 
\mathcal{E}_{\mathcal{B}_{e_{j^1}}^1}^{\hat{\mathcal{B}}}(S_{e_{j^2}^1}(\kappa)) 
\right], \label{eq:S_eta1}\\
&S_{\mathsf{B}_{j^2}\mathsf{C}_{j^2}}(\kappa+1)=\mathcal{P}_
{\mathcal{B}_{\mathsf{B}_{j^2}}\cup \mathcal{B}_{\mathsf{C}_{j^2}}}^{\hat{\mathcal{B}}}\left[
\mathcal{E}_
{\eta_{j^1}}^{\hat{\mathcal{B}}}(S_{\eta_{j^1}}(\kappa))\cap \right. \nonumber\\
&\left. 
\mathcal{E}_
{\mathcal{B}_{\mathsf{B}_{j^2}}\cup \mathcal{B}_{\mathsf{C}_{j^2}}}^{\hat{\mathcal{B}}}(S_{\mathsf{B}_{j^2}\mathsf{C}_{j^2}}(\kappa))\cap 
\mathcal{E}_
{\mathcal{B}_{e_{j^2}^1}}^{\hat{\mathcal{B}}}(S_{e_{j^2}^1}(\kappa)) 
\right], \label{eq:S_mathsf_B_j2_mathsf_C_j2}\\
&S_{e_{j^2}^1}(\kappa+1)=\mathcal{P}_{\mathcal{B}_{e_{j^2}^1}}^{\hat{\mathcal{B}}}\left[
\mathcal{E}_{\mathcal{B}_{\eta_{j^1}}}^{\hat{\mathcal{B}}}(S_{\eta_{j^1}}(\kappa))\cap \right. \nonumber\\
&\left. 
\mathcal{E}_
{\mathcal{B}_{\mathsf{B}_{j^2}}\cup \mathcal{B}_{\mathsf{C}_{j^2}}}^{\hat{\mathcal{B}}}(S_{
\mathsf{B}_{j^2}\mathsf{C}_{j^2}}(\kappa))\cap 
\mathcal{E}_{\mathcal{B}_{e_{j^2}^1}}^{\hat{\mathcal{B}}}(S_{e_{j^2}^1}(\kappa)) 
\right], \label{eq:S_e_j2^1}
\end{align}
\end{subequations}
with $S_{\eta_{j^1}}(0)=\bigcap_{j\in\eta_{j^1}}\mathcal{E}_{\mathcal{B}_j}^{\mathcal{B}_{\eta_{j^1}}}(S_j^*)$, 
\newline
$S_{\mathsf{B}_{j^2}\mathsf{C}_{j^2}}(0)=\bigcap_{j\in (\mathsf{B}_{j^2}\cup \mathsf{C}_{j^2})} \mathcal{E}_{\mathcal{B}_j}^{\mathcal{B}_{\mathsf{B}_{j^2}}\cup \mathcal{B}_{\mathsf{C}_{j^2}}}(S_j^*)$ and $S_{e_{j^2}^1}(0)=S_{e_{j^2}^1}^*$.
Let us consider any  element $\tilde{s}_{j^1}^*\in \tilde{S}_{j^1}^*\subseteq S_{j^1}^*$ (fixed point of  \eqref{eq:distributed_extrusion_generated_set} according to Lemma \ref{lem:fixed_point}).
From the fixed point equations (Lemma \ref{lem:fixed_point} point i.) it follows that there exist vectors $\tilde{s}_{\eta_{j^1}}^*\in S_{\eta_{j^1}}^*$ and $\tilde{s}_{\mathcal{M}_{j^2}}^*\in S_{\mathcal{M}_{j^2}}^*\subseteq \Bbb{R}^{|\mathcal{B}_{\mathcal{M}_{j^2}}|}$, with $\mathcal{B}_{\mathcal{M}_{j^2}}=\bigcup_{j\in\mathcal{M}_{j^2}}\mathcal{B}_j$ according to which the following positions can be considered: $\tilde{s}_{\mathsf{A}_{j^1}}^*:=\mathcal{P}_{\mathcal{B}_{\mathsf{A}_{j^1}}}^{\mathcal{B}_{\eta_{j^1}}}(\tilde{s}_{\eta_{j^1}}^*)$, 
$\tilde{s}_{\mathsf{B}_{j^2}}^*:=\mathcal{P}_{\mathcal{B}_{\mathsf{B}_{j^2}}}^{\mathcal{B}_{\eta_{j^1}}}(\tilde{s}_{\eta_{j^1}}^*)=\mathcal{P}_{\mathcal{B}_{\mathsf{B}_{j^2}}}^{\mathcal{B}_{\mathcal{M}_{j^2}}}(\tilde{s}_{\mathcal{M}_{j^2}}^*)$, 
$\tilde{s}_{\mathsf{D}_{j^1}}^*:=\mathcal{P}_{\mathcal{B}_{\mathsf{D}_{j^1}}}^{\mathcal{B}_{\eta_{j^1}}}(\tilde{s}_{\eta_{j^1}}^*)$,
$\tilde{s}_{\mathsf{C}_{j^2}}^*:=\mathcal{P}_{\mathcal{B}_{\mathsf{C}_{j^2}}}^{\mathcal{B}_{\mathcal{M}_{j^2}}}(\tilde{s}_{\mathcal{M}_{j^2}}^*)$,
 $\tilde{s}_{e_{j^2}^1}^{\prime\prime}:=\mathcal{P}_{\mathcal{B}_{e_{j^2}^1}}^{\mathcal{B}_{\mathcal{M}_{j^2}}}(\tilde{s}_{\mathcal{M}_{j^2}}^*)$. Such a reasoning holds for any point 
 $\tilde{s}_{\eta_{j^1}}^*\in S_{\eta_{j^1}}^*$ and $\tilde{s}_{\mathcal{M}_{j^2}}^*\in S_{\mathcal{M}_{j^2}}^*\subseteq \Bbb{R}^{|\mathcal{B}_{\mathcal{M}_{j^2}}|}$. 
 From Lemma \ref{lem:extracted_set}, point $\tilde{s}_{\eta_{j^1}^*}$ induces on set $S_{e_{j^2}^1}(0)$ a point $\tilde{s}_{e_{j^2}^1}^{\prime}$. It is easy to notice that $\tilde{s}_{e_{j^2}^1}^{\prime}=\tilde{s}_{e_{j^2}^1}^{\prime\prime}$. Indeed, if by contradiction we suppose  $\tilde{s}_{e_{j^2}^1}^{\prime}\neq \tilde{s}_{e_{j^2}^1}^{\prime\prime}$, then $\tilde{s}_{\eta_{j^1}^*}$ would not be, conversely on what stated,  a fixed point. 

In the exact same way we can prove the implication (i)$\Rightarrow$(ii) for $j^1, \dots, j^M\in \eta_{j^1}$, while the implication (ii)$\Rightarrow$(i) is trivial.

However, we still need to prove the first part of the theorem.
According to what have been shown so far, it follows that any nonempty subset $\tilde{S}_{j^1}^*$ induces nonempty subsets $\tilde{S}_j^*$, with $j\in\eta_{j^1}\cup \mathsf{B}_{j^2}\cup\mathsf{C}_{j^2}\cup\{e_{j^2}^1\}$ such that 
\begin{equation}\label{eq:intersection_extrusions1_thm_distributed_extrusion_algorithm_j^1_partial}
\bigcap_{j\in(\eta_{j^1}\cup \mathsf{B}_{j^2}\cup\mathsf{C}_{j^2}\cup\{e_{j^2}^1\})}\mathcal{E}_{\mathcal{B}_j}^{\bar{\mathcal{B}}}(\tilde{S}_j^*)\neq\emptyset. 
\end{equation}

To show that there also exists a nonempty subset $\tilde{S}_{e_{j^2}^2}^*\subseteq S_{e_{j^2}^2}^*$, we consider the update $\mathsf{E}_{j^2}^\prime \leftarrow \mathsf{E}_{j^2}^\prime\cup\{e_{j^2}^1\}$ and proceed as in the derivation of \eqref{eq:eta_j1}$\--$\eqref{eq:mathsf_A_j1} and \eqref{eq:S_eta1}$\--$\eqref{eq:S_e_j2^1}, by replacing $e_{j^2}^1$ with $e_{j^2}^2$. 
Similarly, the same approach is considered for all the elements in $\mathsf{E}_{j^2}$, i.e., $e_{j^2}^h$ with $h=1, \dots, |\mathsf{E}_{j^2}|$, thus extending the nonempty extract set relation from \eqref{eq:intersection_extrusions1_thm_distributed_extrusion_algorithm_j^1_partial}
to 
\begin{equation}
\label{eq:intersection_extrusions1_thm_distributed_extrusion_algorithm_j^1}
\bigcap_{j\in(\eta_{j^1}\cup \mathcal{M}_{j^2})}\mathcal{E}_{\mathcal{B}_j}^{\bar{\mathcal{B}}}(\tilde{S}_j^*)\neq\emptyset. 
\end{equation}

Having completed the iterative procedure for $j^1$ and $j^2$, we highlight that the same steps described so far can be extended to node $j^3$ by considering $\eta_{j^2}=\eta_{j^1}\cup\mathcal{M}_{j^2}$ and shifting by one the indices the equations \eqref{eq:mathsf_E_j2_prime}$\--$\eqref{eq:S_e_j2^1}. Similarly from $j^3$ to $j^4$ and to all the nodes in the sequence.

More formally, by considering the initialization $\eta_{j^0}=\emptyset$, we generalize the procedure described before to all nodes in the sequence $j^2,\dots,j^N$ evaluating for $i=1,\dots,N-1$, the following update equations
\begin{subequations}
\begin{align}
\mathsf{E}_{j^{i+1}}^{\prime}&=\emptyset, \label{eq:mathsf_E_j^i_plus_1_prime_empty}\\
\mathsf{E}_{j^{i+1}}^{\prime}& \leftarrow \mathsf{E}_{j^{i+1}}^{\prime}\cup \{e_{j^{i+1}}^h\}, \,\,\, h=1,\dots, |\mathsf{E}_{j^{i+1}}|,
\end{align}
\end{subequations}
where $\mathsf{E}_{j^{i+1}}=\{e_{j^{i+1}}^1,e_{j^{i+1}}^2,\dots,e_{j^{i+1}}^{|\mathsf{E}_{j^{i+1}}|}\}$, and with the iterative equations
\begin{subequations}
\begin{align}
&\eta_{j^i} :=\eta_{j^{i-1}}\cup\mathcal{M}_{j^i}, \label{eq:eta_j_i}\\
&\mathsf{B}_{j^{i+1}} :=\mathcal{M}_{j^{i+1}}\cap\eta_{j^i}, \label{eq:mathsf_B_j_i_plus_one}\\
&\mathsf{E}_{j^{i+1}} :=\left\{j\in (\mathcal{M}_{j^{i+1}}\backslash\mathsf{B}_{j^{i+1}}): \mathcal{M}_j\cap(\eta_{j^i}\backslash\mathsf{B}_{j^{i+1}})\neq \emptyset \right\}, \label{eq:mathsf_E_j_i_plus_one}\\
&\mathsf{C}_{j^{i+1}}  :=(\mathcal{M}_{j^{i+1}}\backslash(\mathsf{B}_{j^{i+1}}\cup \mathsf{E}_{j^{i+1}}))\cup\mathsf{E}_{j^{i+1}}^\prime , \label{eq:mathsf_C_j_i_plus_one}\\
&\mathsf{D}_{j^i}  :=\left\{j\in (\eta_{j^i}\backslash\mathsf{B}_{j^{i+1}}): \mathcal{M}_j\cap \mathsf{E}_{j^{i+1}}\neq \emptyset \right\},\label{eq:mathsf_D_j_i}\\
&\mathsf{A}_{j^i}  := \eta_{j^i}\backslash(\mathsf{B}_{j^{i+1}}\cup \mathsf{D}_{j^i}). \label{eq:mathsf_A_j_i}
\end{align}
\end{subequations}

The iterative procedure described so far allows to state that, starting from $j^1$, $\tilde{S}_{j^1}^*$ induces a subset $\tilde{S}_{j^2}^*$ such that $\mathcal{E}_{\mathcal{B}_{j^1}}^{\bar{\mathcal{B}}}(\tilde{S}_{j^1}^*)\cap \mathcal{E}_{\mathcal{B}_{j^2}}^{\bar{\mathcal{B}}}(\tilde{S}_{j^2}^*)\neq \emptyset$. 
For the sake of simple notation, let $\tilde{Y}_{j^1}^*:=\mathcal{E}_{\mathcal{B}_{j^1}}^{\bar{\mathcal{B}}}(\tilde{S}_{j^1}^*)$ and 
 $\tilde{Y}_{j^2}^{*,1}:=\mathcal{E}_{\mathcal{B}_{j^2}}^{\bar{\mathcal{B}}}(\tilde{S}_{j^2}^*)$, where the superscript $1$ recalls that $\tilde{S}_{j^2}^*$ has been induced by $\tilde{S}_{j^1}^*$. 
Then, both $\tilde{S}_{j^1}^*$ and $\tilde{S}_{j^2}^*$ induce the nonempty set $\tilde{S}_{j^3}^*$ such that  $\tilde{Y}_{j^1}^*\cap\tilde{Y}_{j^2}^{*,1}\cap\tilde{Y}_{j^3}^{*,1,2} \neq \emptyset$, with $\tilde{Y}_{j^3}^{*,1,2}$ defined analogously to  $\tilde{Y}_{j^2}^{*,1}$, i.e., $\tilde{Y}_{j^3}^{*,1,2}:=\mathcal{E}_{\mathcal{B}_{j^3}}^{\bar{\mathcal{B}}}(\tilde{S}_{j^3}^*)$.
Iterating along the path implies that, for any $2\leq i\leq N$, we have
\begin{equation}\label{eq:path_iteration_extracted_extrusion}
\tilde{Y}_{j^1}^*\cap\tilde{Y}_{j^2}^{*,1}\cap\dots, \cap \tilde{Y}_{j^i}^{*,1,2, \dots, i-1} \neq \emptyset.
\end{equation}
This concludes the proof of the first part of the theorem.
\end{pf}

\begin{pf}\textbf{(Theorem \ref{thm:thm_distributed_extrusion_algorithm})}
By combining Lemma \ref{lem:succession_nested_inclusions} (set sequence) and Lemma \ref{lem:fixed_point}, the theorem will be proven by showing that $S_i^*=\bar{S}_i, \,\, \forall i=1, \dots, N$. We first focus on the case of connected graph $\mathcal{G}_c=(R,E)$ and then extended the proof to the case of non connected graphs.

Taking into account \eqref{eq:S_i_star} and Lemma \ref{lem:succession_nested_inclusions}, we immediately have that $\bar{S}_i\subseteq S_i^*$.
 
To prove the equality between the two sets let us suppose, by  contradiction, that at least for one $i$, we have $S_i^* \neq \bar{S}_i$ and so that the strict inclusion $\bar{S}_i\subset S_i^*$ holds. From such strict inclusion assumption,  
by setting 
\begin{equation}\label{eq:S_tilde_i_star_difference}
\tilde{S}_i^*=S_i^*\backslash\bar{S}_i,
\end{equation}
it follows that $\tilde{S}_i^*\neq \emptyset$.

For the sake of brevity, we use the same simplified notation adopted in the proof of Theorem \ref{thm:iterative_extracted_sets}. Specifically, we set $Y_j^*=\mathcal{E}_{\mathcal{B}_j}^{\bar{\mathcal{B}}}(S_j^*)$ and $\tilde{Y}_j^*=\mathcal{E}_{\mathcal{B}_j}^{\bar{\mathcal{B}}}(\tilde{S}_j^*)$ (the latter with the possibility of reporting in the superscript the inducing subsets, as done for \eqref{eq:path_iteration_extracted_extrusion}). 
The subsets $\tilde{S}_j^*$, with $j\neq i$, will be determined next. 

The following relation holds
\begin{equation}\label{eq:extraxted_extrusion_general_centralized_intersection}
\tilde{Y}_i^*\cap \bigcap_{j=1, j\neq i}^N Y_j^*=\emptyset.
\end{equation}
Indeed, we have 
\begin{align*}
\tilde{Y}_i^* &\cap \bigcap_{j=1, j\neq i}^N Y_j^*=\\
\tilde{Y}_i^* & \cap \bigcap_{j=1}^N Y_j^*=\\
\tilde{Y}_i^* & \cap \bar{S}=\emptyset ,
\end{align*}
where, in the last step we applied Lemma \ref{lem:fixed_point} point ii. and considered that, in case of nonempty intersection, we would have had points $\bar{s}_i\in\mathcal{P}_{\mathcal{B}_i}^{\bar{\mathcal{B}}}(\bar{S})$ such that $\bar{s}_i\in\tilde{S}_i^*$ and $\bar{s}_i\in\bar{S}_i$, contradicting the \eqref{eq:S_tilde_i_star_difference}.   

Let us now order the nodes as $j^1,j^2,\dots, j^N$ according to a connected path and with $j^1=i$. We call
\begin{equation*}
\mathcal{A}_0=\tilde{Y}_{j^1}^*\cap Y_{j^2}^*\cap Y_{j^3}^*\cap \dots \cap Y_{j^N}^*.
\end{equation*}

We rewrite \eqref{eq:extraxted_extrusion_general_centralized_intersection} as $\mathcal{A}_0=\emptyset$. Since $\tilde{S}_{j^1}^*\neq \emptyset$, it induces a nonempty subset to $S_{j^2}^*$ along the connected path, according Theorem \ref{thm:iterative_extracted_sets}. 
Let now consider
\begin{equation*}
\mathcal{A}_1=\tilde{Y}_{j^1}^*\cap \tilde{Y}_{j^2}^{*,1}\cap Y_{j^3}^*\cap \dots \cap Y_{j^N}^*.
\end{equation*}

We obviously have $\mathcal{A}_1\subseteq \mathcal{A}_0$. Considering now the induced $\tilde{Y}_{j^3}^{*,1,2}$ and defining $\mathcal{A}_2$ in analogous way, we have $\mathcal{A}_2\subseteq \mathcal{A}_1\subseteq \mathcal{A}_0$.
Generalizing the iterations up to node $j^{i+1}$, with $i=1, \dots, N-1$, and by iteratively applying Theorem \ref{thm:iterative_extracted_sets} we have
\begin{equation*}
\mathcal{A}_i=\tilde{Y}_{j^1}^*\cap \tilde{Y}_{j^2}^{*,1}\cap \tilde{Y}_{j^3}^{*,1,2}\cap \dots \cap
\tilde{Y}_{j^{i+1}}^{*,1,2,\dots,i}\cap
Y_{j^{i+1}}^*\cap\dots\cap
Y_{j^{N}}^*,
\end{equation*}
with $\mathcal{A}_i\subseteq \mathcal{A}_{i-1}$.

Iterating the process up to node $j^{N}$ we have that, according to \eqref{eq:path_iteration_extracted_extrusion}, $\mathcal{A}_{N-1}\neq \emptyset$. 
On the other hand, since from \eqref{eq:extraxted_extrusion_general_centralized_intersection} $\mathcal{A}_0=\emptyset$, we have 
\begin{equation}\label{eq:contraddiction_A0}
\emptyset\neq \mathcal{A}_{N-1}\subseteq \mathcal{A}_0=\emptyset.
\end{equation}
The contradiction implies that $\bar{S}_i=S_i^*$.

The poof can be easily extended to the case of non connected graphs. Indeed suppose that, without loss of generality, the graph is partitioned in $p$ connected components, i.e., $\mathcal{G}_{c1}=(R_1,E_1), \mathcal{G}_{c2}=(R_2,E_2), \dots, \mathcal{G}_{cp}=(R_p,E_p)$, with $N_1, N_2, \dots, N_p$ elements each. 

Notice that, for any $\mathcal{G}_{c_m}=(R_m, E_m)$ and $\mathcal{G}_{c_l}=(R_l, E_l)$ different components, i.e., with $m\neq l$, we always have $\tilde{Y}_{j}^*\cap \tilde{Y}_{h}^*\neq \emptyset$ with $j\in R_m$ and $h\in R_l$ since $\mathcal{B}_j\cap \mathcal{B}_h=\emptyset$.
Hence, for the whole graph we can manipulate each connected component independently and, at the same time, write again $\mathcal{A}_0=\emptyset$.  
By using similar arguments as for \eqref{eq:contraddiction_A0}, we obtain
\begin{equation*}
\emptyset \neq 
\mathcal{A}_{N_1-1}\cap
\mathcal{A}_{N_2-1}\cap\dots\cap 
\mathcal{A}_{N_p-1}\subseteq \mathcal{A}_0=\emptyset.
\end{equation*}

The above contradiction implies $\bar{S}_i=S_i^*$.
\end{pf}

\begin{pf}\textbf{(Lemma \ref{lem:reachability_overall_system})}
The proof is trivial. It suffices to notice that the solution of system \eqref{eq:system_of_constraints_overall} satisfies the conditions in Definition \ref{def:reachability_H_steps} and vice-versa.
\end{pf}

\begin{pf}\textbf{(Lemma \ref{lem:pre_overall_system})}
Consider the set $\bar{S}_k:=\mathcal{P}_{\bar{\mathcal{B}}_{x,0}}^{\bar{\mathcal{B}}^H}(\bar{S}_{kh})$. It is trivial to notice that for all the $\bar{\mathbf{z}}\in \bar{S}_k$ there exist $\mathbf{z}\in \bar{S}_{kh}$ such that $\mathcal{P}_{\bar{\mathcal{B}}_{x,0}}^{\bar{\mathcal{B}}^H }(\mathbf{z})\in \bar{S}_k$. Since $\bar{S}_{kh}$ is the solution of \eqref{eq:system_of_constraints_overall}, then $\bar{S}_k  \rightrsquigarrow \bar{S}_h$.
 
On the other hand, by setting $\breve{S}_k=\Bbb{R}^n\backslash\bar{S}_k$, we trivially have $\breve{S}_k  \nrightrsquigarrow \bar{S}_h$. Indeed, $\forall \mathbf{z}\in \Bbb{R}^{(H+1)(n+m)}$ such that $\mathcal{P}_{\bar{\mathcal{B}}_{x,0}}^{\bar{\mathcal{B}}^H}(\mathbf{z})\in \breve{S}_k$, we have that $\mathbf{z}\notin \bar{S}_{kh}$. Therefore, $\bar{S}_k=\mathrm{Pre}(\bar{S}_h)$. 
\end{pf}

\begin{pf}\textbf{(Lemma \ref{lem:bar_S_kh_extrusion_generated})}
The proof is trivial since the left hand side of equation \eqref{eq:bar_S_kh_extrusion_generated} is the solution of system \eqref{eq:system_of_constraints_overall}, while the right hand side is the solution of the system obtained by considering the \eqref{eq:system_of_constraints_i} for all $i=1,\dots, N$. The latter is an equivalent way of expressing system \eqref{eq:system_of_constraints_overall}.
\end{pf}

\begin{pf}\textbf{(Theorem \ref{thm:distributed_computation_bar_S_kh})}
The proof is a direct result of the application of Theorem \ref{thm:thm_distributed_extrusion_algorithm}.
\end{pf}

\begin{pf}\textbf{(Theorem \ref{thm:distributed_Pre})}
It is immediate to notice that, for all $i=1,\dots, N$,  $\bar{S}_{k,i}=
\mathcal{P}_{\mathcal{B}_{x,0,i}}^{\bar{\mathcal{B}}^H}(\bar{S}_{kh})$. Indeed, $\bar{S}_{k,i}=\mathcal{P}_{\mathcal{B}_{x,0,i}}^{\bar{\mathcal{B}}_{x,0}}(\bar{S}_k)=
\mathcal{P}_{\mathcal{B}_{x,0,i}}^{\bar{\mathcal{B}}_{x,0}}\left[
\mathcal{P}_{\bar{\mathcal{B}}_{x,0}}^{\bar{\mathcal{B}}^H}(\bar{S}_{kh})
\right]
=
\mathcal{P}_{\mathcal{B}_{x,0,i}}^{\bar{\mathcal{B}}^H}(\bar{S}_{kh})$.

We will prove that $\underline{S}_{k,i}=\bar{S}_{k,i}$ via showing the double inclusion $\bar{S}_{k,i}\subseteq \underline{S}_{k,i}$ and $\underline{S}_{k,i}\subseteq\bar{S}_{k,i}$, with $i=1, \dots, N$.

Let us call $\underline{S}_{k}=\bigcap_{j=1}^N\mathcal{E}_{\mathcal{B}_{x,0,j}}^{\bar{\mathcal{B}}_{x,0}}(\underline{S}_{k,j})$. 
First of all, notice that by construction $\bar{S}_k\subseteq\underline{S}_k$. Since for every point $\bar{s}_{k,i}$ there exist points $\bar{s}_{k}\in\bar{S}_{k}$ such that $\mathcal{P}_{\mathcal{B}_{x,0,i}}^{\bar{\mathcal{B}}_{x,0}}(\bar{s}_k)=\bar{s}_{k,i}$ and since $\bar{s}_{k}\in\underline{S}_{k}$, then $\bar{s}_{k,i}\in \underline{S}_{k,i}$. Therefore, $\bar{S}_{k,i}\subseteq\underline{S}_{k,i}$. 

To prove that for every point $\underline{s}_{k,i}\in\underline{S}_{k,i}$ we also have $\underline{s}_{k,i}\in\bar{S}_{k,i}$ it suffices to notice that there exist points $\bar{s}_{kh}\in\bar{S}_{kh}$ such that $\mathcal{P}_{\mathcal{B}_{x,0,i}}^{\bar{\mathcal{B}}^H}(\bar{s}_{kh})=\underline{s}_{k,i}$. Since the projection of such points $\bar{s}_{kh}$ onto the space generated by the axis set $\bar{\mathcal{B}}_{x,0}$ gives points $\bar{s}_k$ such that $\mathcal{P}_{\mathcal{B}_{x,0,i}}^{\bar{\mathcal{B}}_{x,0}}(\bar{s}_k)=\underline{s}_{k,i}$, then $\underline{s}_{k,i}\in \bar{S}_{k,i}$. This concludes the proof. 
\end{pf}

\begin{pf}\textbf{(Theorem \ref{thm:distributed_flow})}
The proof follows analogous steps as for the proof of Theorem \ref{thm:distributed_Pre}.
\end{pf}

\end{document}